\documentclass[11pt]{article}
\usepackage{amsfonts}
\usepackage{indentfirst}
\usepackage{graphics}
\usepackage{amssymb}
\usepackage{color,xcolor}
\usepackage{graphicx,epstopdf}
\usepackage{epsfig}
\usepackage{subfigure}
\usepackage{caption}
\usepackage{mathrsfs}
\usepackage{bbm}
\usepackage{chngpage}
\usepackage{cite}
\usepackage{multirow}
\usepackage{multicol}
\usepackage{dsfont}
\usepackage{bm}
\usepackage{amssymb,amsmath,amsthm,mathrsfs}
\usepackage{color}
\usepackage{booktabs}
\usepackage{algpseudocode}
\usepackage{algcompatible}
\usepackage{float}
\usepackage{xspace}
\usepackage{algpseudocode}
\usepackage{pifont}

\newcommand{\ben}{\begin{eqnarray}}
\newcommand{\een}{\end{eqnarray}}
\newcommand{\bea}{\begin{array}}
\newcommand{\eea}{\end{array}}
\allowdisplaybreaks

\newtheorem{theorem}{Theorem}[section]

\newtheorem{algorithm}[theorem]{Algorithm}
\theoremstyle{definition}

\theoremstyle{remark}
\newtheorem{remark}[theorem]{Remark}
\newtheorem{lemma}{Lemma}[section]

\newcommand{\p}{\partial}

\def\p{\partial}

\newcommand{\xx}{{\bf x}}

\newcommand{\myvec}{\stackrel{\raisebox{-2pt}[0pt][0pt]{\small$\rightharpoonup$}}}

\begin{document}
	\title	{Accelerated numerical algorithms for steady states of Gross-Pitaevskii equations coupled with microwaves}
	
	\author{Di Wang\footnote{Beijing Computational Science Research Center, Beijing 100193, China. Email: di\_wang@csrc.ac.cn.} \ and \ Qi Wang\footnote{Department of Mathematics, University of South Carolina, Columbia, SC 29208, USA. Email: qwang@math.sc.edu.}}

\date{}	
\maketitle
\textbf{Abstract.} We present two accelerated numerical algorithms for  single-component and binary Gross-Pitaevskii (GP) equations coupled with microwaves (electromagnetic fields) in steady state. One is  based on a normalized gradient flow formulation, called the ASGF method, while the other on a perturbed,  projected conjugate gradient approach for the nonlinear constrained optimization, called the PPNCG method.   The coupled GP equations are nonlocal in space, describing pseudo-spinor Bose-Einstein condensates (BECs) interacting with an electromagnetic field. Our interest in this study is to develop efficient, iterative numerical methods for steady symmetric and central vortex states of the nonlocal GP equation systems. In the algorithms, the GP equations are discretized by a Legendre-Galerkin spectral method in a polar coordinate in two-dimensional (2D) space.  The new algorithms are shown to outperform the existing ones through a host of benchmark examples, among which the PPNCG method performs the best.  Additional numerical simulations of the central vortex states are provided to  demonstrate the usefulness and efficiency  of the new algorithms.
\\


\section{Introduction}

The centrepiece of studies on BECs lies in the study of  quantized vortices, which are  building blocks of quantum turbulence \cite{PhysRevX.5.021015, PhysRevLett.110.215302, HadzibabicBerezinskii, PhysRevLett.96.020404, PhysRevA.91.013612, PhysRevA.91.053615, PhysRevLett.113.165302, PhysRevLett.111.235301, PhysRevLett.103.045301}. In addition to creating traps and optical lattices \cite{ADAMS1994143, RevModPhys.80.1215, GRYNBERG2001335}, various optical patterns associated with quantum vortices have potential applications in the field of quantum data processing \cite{PhysRevLett.98.203601, PhysRevLett.97.170406}. In this study, we explore accelerated numerical algorithms for computing 2D steady vortices in a binary atomic BEC interacting with a (electromagnetic) microwave field.

BECs at temperature T much lower than the critical condensation temperature $T_c$ are usually well modelled by a nonlinear Schr\"odinger equation (NLSE) for the macroscopic wave function known as the Gross-Pitaevskii (GP)  equation \cite{bao2018mathematical, BaoMathematical, RevModPhys.71.463, KAWAGUCHI2012253}. One of the fundamental issues in the study of the equation is to study the equation's steady states of certain properties, for instance the ground state and excited states. The ground state is usually defined as the minimizer of the energy functional  under the normalization constraint for the wave function.

The steady state solution whose corresponding energy is larger than that of the ground state is usually called an excited state. Among the exited states, there are some vortex steady states of winding number (or topological charge) $S>0$ (which will be defined precisely in the text).  In BECs in a rotational frame with an angular velocity, self-trapped vortex annuli (VA) with large values of winding number S (giant VA) not only are a subject of fundamental interest in quantum physics, but are also sought for various applications, such as quantum information processing and storage \cite{PhysRevLett.100.223601, PhysRevLett.98.203601, PhysRevLett.97.170406}. To study these states and their properties, an important prerequisite is to find an efficient and accurate solver for the central vortex states, that is, the first ground state of the corresponding Hamiltonian with the vortex in the rotational center \cite{Madison2000, Rokhsar1997}.

In the last two decades, there have been a plethora of numerical methods developed to compute ground states of BECs, including normalized gradient flow methods based on the Hamiltonian (the energy) of the GP equation  \cite{baoComputing, BAO2006836, PhysRevE.62.7438, MR2366691, WANG2014473, Antoine2015, ZENG2009854, ANTOINE2014509, ADHIKARI200091, Baye2010, Cerimele2000, MR4214366, MR4216953, ZHUANG201972}, and methods for the nonlinear eigenvalue problem (see e.g., \cite{MR2679792, CHEN20112222, MR3819575, wang_jeng_chien_2013, DION2007787, MR3946658} and references therein) based on time-independent GP equation as well as constrained optimization techniques \cite{CALIARI2009349, BAO2003230, DANAILA20106946, MR2684722, ANTOINE201792, MR4129007, MR3704854, MR3731036}. The normalized gradient flow strategy is considered from the PDE perspective, leading to numerical algorithms for a dissipative system. Among these methods, the gradient flow with discrete normalization (GFDN) method (also known as the imaginary time evolution method) \cite{baoComputing, BAO2006836, PhysRevE.62.7438}, the continuous normalized gradient flow (CNGF) method \cite{MR2366691, baoComputing, WANG2014473} are two main approaches. Some error estimates \cite{Faou2018} and numerical observations \cite{BAO2006836} noted that GFDN method can converge to a spurious ground state solution with errors depending on the time step size. As an improvement, the GFDN method with imposed explicit Lagrange multiplier terms (GFLM) \cite{MR4214366}, which can be viewed as a special temporal discretization for the CNGF method, is proposed to mitigate the situation. The constrained optimization approaches for the nonlinear eigenvalue problem  include the finite element method directly minimizing the energy functional \cite{BAO2003230}, the Sobolev gradient method \cite{MR2684722}, the regularized Newton method \cite{MR3704854}, the Riemannian optimization method \cite{MR3731036, MR4129007}, the preconditioned, nonlinear conjugate gradient (PNCG) method \cite{ANTOINE201792}, and so on.  For computing symmetric and central vortex states in BECs, a generalized-Laguerre–Hermite pseudospectral method without  truncating the computational domain \cite{BAO20089778} is also proposed.  In \cite{MR2132826}, symmetric and central vortex states in rotating BECs are numerically investigated.

\begin{figure}
	\centering
	\includegraphics[width=12cm]{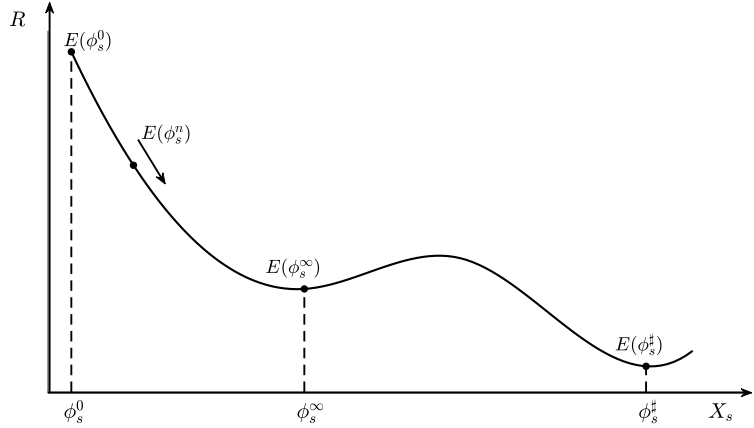}
	\label{fig}
	\caption{The accelerated momentum-based strategy.}
	\label{figure}
\end{figure}

Notice that most Hamiltonians in interacting boson systems like BECs are quite involved, in which the energy functional is non-convex so that the energy landscape presents multiple local shallower minima. In this case, the global uniqueness of the ground state or central vortex state solution is very difficult  to obtain numerically. Under the circumstance, the gradient flow strategy, which is essentially the steepest descent method, may be inadequate. It is known that gradient flows have optimal worst-case complexity for convergence to stationary points, but are strongly attracted to local minima. To mitigate this, one resorts to the accelerated momentum-based methods by adding the inertia back to the gradient flow as follows:
\begin{align} \label{accelerated momentum-based methods}
\left\{\begin{array}{l}
\ddot{\phi}_S(t)+\alpha_1(t) \dot{\phi}_S(t)+\alpha_2(t) \nabla E(\phi_S(t))=0, \\
\phi_S(0)=\phi_{S}^{0}, \quad \dot{\phi}_S(0)=\dot{\phi}_{S}^{0}.
\end{array}\right.
\end{align}
The Polyak’s heavy ball method \cite{MR1753136, Bostan2018, Xu2018} and Nesterov’s accelerated gradient descent method \cite{Nesterov2018, Nesterov1983, Betancourt2018, Wibisono2016} are two good examples. With the addition of inertia, the total energy of the system is augmented with the so-called kinetic energy so that  one can derive a "global" strategy for the numerical computation of the local minima of $E(\phi_S)$ (See Fig. \ref{figure}). In this way, the total energy exhibits global decay (non-oscillatory) in time while the Hamiltonian of GP equation may exhibit under-damped oscillations around the equilibrium, creating opportunities to hop over low barriers to reach lower energy levels.
The steady state solution of the augmented system is completely determined by its initial position and velocity. Then, one can play with
the initial velocity $\dot{\phi}_{S}^{0}$ to reach asymptotically different critical points (local minima). Decades of empirical experience suggests that momentum methods are capable of exploring multiple local minima, which gives them advantages over purely dissipative gradient flows. Moreover, recent theoretical results have demonstrated that momentum methods can escape saddle points faster than standard gradient descent methods \cite{jin2018, ONeill2017}, providing further evidence of their value in nonconvex optimization.

The gradient descent approach is rooted in the dissipative  PDE theory, where the energy functional as the free energy of the gradient flow model decays along a "smooth" path or trajectory.  In practical numerical implementations, for various treatments of   $\nabla E(\Phi_S(t))$ in Eq. \eqref{accelerated momentum-based methods}, $\nabla E(\Phi_S(t))$ may depend on state $\Phi_S(t)$ approximated at different time levels. As the result,  the gradient-flow equations  discretized with respect to pseudotime t may no longer be in the form of a discrete gradient flow of the energy functional. So, the resulting methods typically do not preserve the gradient-flow structure at the discrete level.  This indicates that the structure-preserving strategy for the gradient flow discretization is not important for developing iterative numerical schemes for the minimization. Another potential drawback of such an approach is that solutions of the gradient flow model are in general critical points of the free energy, but are not necessarily minima (i.e., they can be saddle points). From this observation, we notice that preserving the gradient flow structure is secondary in designing an iterative optimization algorithm for the free energy. Several algorithms for the constrained optimization on Riemannian manifolds have been developed, in which constrained analogues of gradient, conjugate gradient and Newton’s algorithms are derived \cite{Absil2009, Edelma1998}.  These provide alternative approaches for us to follow in this study. Motivated by these developments, we will employ a perturbed,  preconditioned, nonlinear conjugate gradient method (PPNCG) on the manifold $\mathcal{S}$ that guarantees the $L_2$ norm constraint for the solution of the nonlocal GP equation.

We note that vortex steady states of the coupled GP model often exhibit  fine spatial structures, imposing strong requirements on the spatial discretization of the PDE system. To retain the required spatial resolution near the fine structures, we adopt the highly accurate Legendre-Galerkin spectral method \cite{Shen1997} to discretize the PDE system in space. In time, we develop two strategies, one is in the accelerated momentum method, called ASGF method, and the other in the projected conjugate gradient method, called PPNCG method. The ASGF strategy hinges on a stabilizer corresponding to a nonlocal inertia "regularization" of the over-damped normalized gradient flow model. This momentum-based method mitigates the strong local-minima attractive nature of the over-damped  gradient flow to facilitate convergence to global minima. In the PPNCG method, we implement a perturbation strategy in the projected conjugate gradient method to avoid saddle points during minimization of the nonconvex Hamiltonian effectively.  The algorithms resulted from both approaches are compared with the existing GFLM method extensively. The numerical results show that new methods perform better than the GFLM method while the PPNCG method outperforms the ASGF method, providing two efficient numerical solvers for solving the steady, coupled, nonlocal GP equations. It is worth noting that the new algorithms not only work well for the nonlocal GP systems, but also efficient for the simpler local GP equations.

In addition to the algorithms we present in this paper, we also devised other algorithms based on several selected high order time discretizations of the spatially, semi-discretized gradient flow systems. These algorithms include algorithms derived from applying second order time discretization, explicit and implicit 4th order Runge-Kutta discretization. None of the resulting iterative schemes outperforms the ASGF algorithm we present in this paper. This indicates that higher order temporal schemes applied to the normalized gradient flow  does not necessarily yield better iterative schemes for steady states of the GP equations.

The rest of this paper is organized as follows. In \S2, we present two new algorithms for a simplified GP model in the case of single-component BECs, detailing the temporal and spatial discretization strategies, and compare the new methods with the existing GFLM method.   In \S 3, we extend the methods to the coupled GP model for binary BECs interacting with microwaves and compare their performance when computing symmetric  and central vortex states.
Finally, we draw the conclusion in \S 4.

\section{Numerical methods for the single-component nonlocal GP equation}

We consider the dimensionless 2D self-trapped single-component nonlocal GP equation in the weak microwave detuning limit \cite{QinStable, WANG2021132852} without the external potential:
\begin{align}\label{GPE}
i \frac{\partial \psi}{\partial t}=\left[-\frac{1}{2} \triangle-\beta |\psi|^{2} -\breve{H}|\psi|^{2}\right] \psi, \quad \xx \in \mathbb{R}^2,
\end{align}
subject to constraint
\begin{align}
\int_{\mathbb{R}^2}|\psi(\xx)|^2d\xx=1,
\end{align}
where  the magnetic field $\breve{H}$ satisfies the following Poisson equation:
\begin{align}
-\triangle \breve{H}= \gamma |\psi|^{2} \quad \xx \in \mathbb{R}^2.\label{eq-H}
\end{align}
We identify the Hamiltonian of the conservative system as follows:
\begin{align}\label{energy functional2}
E(\psi)= \int_{\mathbb{R}^2}\bigg[\frac12|\nabla\psi|^2-\frac{\beta}2|\psi(\xx)|^4-\frac{\breve{H}}2|\psi(\xx)|^4\bigg]d\xx.
\end{align}

To find the symmetric and central vortex steady state solution of \eqref{energy functional2},   we seek the following solution ansatz in polar coordinate $(r,\theta)$:
\begin{align}\label{stationary single}
\psi(\xx)=e^{-i S \theta} \phi_{S}(r),
\end{align}
where $S\in \mathbb{Z}$ is called the winding number and $\phi_S(r)$ is a real-valued function of $r=\sqrt{x^2+y^2}$.  Since the laplace operator is rotational invariant and function $|\psi(\xx)|$ is radially symmetric, it follows from \eqref{eq-H} that magnetic field $\breve{H}$ is radially symmetric, its governing Poisson equation reduces to
\begin{align}\label{poisson R+}
-\frac1{r} \frac{d}{dr}\big(r\frac{d H}{dr}\big)=\gamma |\phi_S|^2, \quad r\in [0,\infty),
\end{align}
in the polar coordinate $(r, \theta)$.
When $S=0$, the solution is call a symmetric state; while $S>0$, it is called a central vortex state.
For this type of solutions, Hamiltonian \eqref{energy functional2} reduces to the following  functional  parameterized by winding number S:
\begin{align}\label{single energy functional_polar}
E(\phi_S)=\pi\int_0^{\infty}\bigg[\left(\left(\phi_{S}^{\prime}\right)^{2} +\frac{S^2}{r^2} \phi_{S}^{2}\right)  -\beta\phi_{S}^4-H \phi_{S}^2\bigg]\ rdr,
\end{align}
subject to  constraint
\begin{align}\label{single polar constraint}
2\pi\int_0^{\infty}\phi_S^2(r)\ r dr=1.
\end{align}

Our objective in this study is to solve for $\Phi_S(r)$ from the nonlocal GP equation. We consider solution $\phi_S(r)$  in the following function space:
\begin{align}
X_S(\mathbb{R^+}): = H^1_S(\mathbb{R^+}) \cap L^4_S(\mathbb{R^+}) \cap L^2_{log,S}(\mathbb{R^+}),
\end{align}
whose various norms are defined as follows:
$$ \|u\|_{L^m_S(\mathbb{R^+})}:=\big(2\pi\int_0^{+\infty}|u(r)|^mr dr\big)^{1/m},\quad m>0,$$
$$ |u|_{H^1_S(\mathbb{R^+})}:=\bigg(2\pi\int_0^{+\infty}\big(|u'(r)|^2+\frac{S^2}{r^2}|u(r)|^2\big)r dr\bigg)^{1/2}, $$
$$\|u\|_{H^1_S(\mathbb{R^+})}:=\bigg(2\pi\int_0^{+\infty}\big(|u'(r)|^2+(\frac{S^2}{r^2}+1)|u(r)|^2\big)r dr\bigg)^{1/2}, $$
and
$$ \|u\|_{L^2_{log,S}(\mathbb{R^+})}:=\big(2\pi\int_0^{+\infty} \ln(1+r) |u(r)|^2r dr\big)^{1/2}.$$

For a given $S$, we denote the symmetric state $\phi_S(r)$ at $S=0$ as $\phi_S^s$ and central vortex state when $S>0$ as  $\phi_S^c$, respectively, both of which minimize $E(\phi_S)$ at respective values of $S$ confined to manifold
\begin{align} \label{single manifold}
\mathcal{S}_1 :=\left\{\phi_S \ \big|\ \|\phi_S\|_{L^2(\mathbb{R}^+)}^{2}=2\pi\int_0^{\infty}\phi_S^2(r)\ r dr=1, E(\phi_S)<\infty\right\}.
\end{align}
The Euler-Lagrange equation when minimizing \eqref{single energy functional_polar} over \eqref{single polar constraint} is given by
\begin{align}\label{single lagrange equation}
\mu\phi_{S}&=\frac12 \frac{\delta E(\phi_S)}{\delta \phi_S} =-\frac12\triangle_{r,S}\phi_{S} -\beta\phi_{S}^2\phi_{S}-H\phi_{S},
\end{align}
where
\begin{align}
\triangle_{r,S} \triangleq \p_r^2+\frac1{r}\p_r-\frac{S^2}{r^2},
\end{align}
and $\mu$  serves as a Lagrange multiplier or nonlinear eigenvalue, which is given by constant
\begin{align}\label{single chemical polar}
\mu(\phi_S)=2\pi\int_0^{\infty}\bigg[\frac12\left(\left(\phi_{S}^{\prime}\right)^{2} +\frac{S^2}{r^2} \phi_{S}^{2}\right) -\beta\phi_{S}^4-H \phi_{S}^2\bigg]\ rdr,
\end{align}
following \eqref{single polar constraint}.

\begin{lemma} (see \cite{WANG2021132852})
	There exists a symmetric state $(S=0)$ and a central vortex state $(S>0)$ of \eqref{single energy functional_polar}  when $\beta< \beta_b^1$, where $\beta_b^1$ is listed in Table \ref{Tab1}.
	\begin{table}
	\centering
	\begin{tabular}{p{0.5cm}p{1cm}p{1cm}p{1cm}p{1cm}p{1cm}p{1cm}p{1cm}p{1cm}p{1cm}{c}}
	\hline
	$\mathrm{S}$ & 0 & 1 & 2 & 3 & 4 & 5 & 6 & 7 \\
	$\beta_{b}^1$ & 5.85 & 24.16 & 44.88 & 66.21 & 87.75 & 109.38 & 131.06 & 152.76 \\
	\hline
	 $\mathrm{S}$ & 8 & 9 & 10 & 11 & 12 & 13 & 14 & 15 \\
	$\beta_{b}^1$ & 174.47 & 196.20 & 217.94 & 239.68 & 261.42 & 283.17 & 304.92 & 326.67 \\
	\hline
\end{tabular}
	\caption{$\beta_b^1$ vs $S$.}
    \label{Tab1}
	\end{table}
	When $\beta> \beta_b^1$, there does not exist any symmetric or central vortex state.
\end{lemma}

By Agmon’s Theorem (see \cite{agmon2014lectures}), it is easy to deduce that
\begin{align}
\phi_S(r) = o(e^{-\alpha r}) \quad \text{as $r \rightarrow \infty$ for every $\alpha > 0$.}
\end{align}
 Hence it's reasonable to truncate the full domain $\mathbb{R}^+$ to the finite domain $U \triangleq [0, R]$ when solving the Euler-Lagrange equation numerically, where $R$ is a sufficiently large positive number. The boundary condition of the solution is given by $\phi_S(R)=0$, together with either $\phi_S(0)=0$ for $S>0$ or  $\frac{d }{dr} \phi_S(0)=0$ for $S=0$.

 \begin{remark}
 	We note that pole condition $\frac{d }{dr} \phi_S(0)=0$ for $S=0$, derived from the parity argument is, however, not part of the essential pole condition for \eqref{single lagrange equation} \cite{Canut1987, Fornberg1995, Shen1997}. Although in most cases there is no harm to impose this extra pole condition, we choose not to do so in our spectral representation since its implementation is more complicated and it may fail to give accurate results in some extreme (but still legitimate) cases.
 \end{remark}

In this case, Eq. \eqref{poisson R+} for the magnetic field $H$ reduces to
\begin{align}\label{poisson}
-\frac1{r} \frac{d}{dr}\big(r\frac{d H}{dr}\big)=\gamma |\phi_S|^2, \quad r\in [0,R],
\end{align}
subject to a Robin boundary condition at $r=R$ \cite{Mauser, WANG2021132852}:
\begin{align}
\frac{d H}{dr}\bigg|_{r=R}=\frac{H(R)}{R ln(R)}.
\end{align}

Next, we present the first numerical method for solving the constrained  minimization problem, called the accelerated, stabilizer-based normalized gradient flow method with Lagrange multipliers (ASGF).

\subsection{Accelerated, stabilized normalized gradient flow (ASGF) method }

We treat the minimization problem for the steady state over  manifold $\mathcal{S}_1$ as a steady state solution of a gradient flow with the Hamiltonian as the free energy of the relaxation dynamics defined in the manifold, i.e.,
\begin{align}
\partial_t \phi_S=-\frac{1}{2}\frac{\delta E}{\delta \phi_S}, \quad \phi_S \in \mathcal{S}_1,
\end{align}
where $\phi_S(\xx,t)$ is treated as a pseudo-time (t) dependent function.
We divide time interval $[0,\infty)$, using  time step $\tau > 0$, into $[t_n, t_{n+1}]$, where $t_n = n\tau$ for $n = 0,1,\cdots, \infty$.
 To deal with the confinement in the manifold, a simple projection step is implemented at the end of each interval.  This method is known as the gradient flow with discrete normalization (GFDN) method \cite{baoComputing}, in which the corresponding PDE and the end-point projection are given  as follows
\begin{align} \label{GFDN}
\begin{split}
&\p_t \phi_S = -\frac12 \frac{\delta E(\phi_S)}{\delta \phi_S} = \frac12\triangle_{r,S}\phi_{S}+\big(\beta\phi_{S}^2+H\big)\phi_{S}, \quad t_{n}<t<t_{n+1}, \quad n \geq 0, \\
&\phi_S\left(r, t_{n+1}\right) \triangleq \phi_S\left(r, t_{n+1}^{+}\right)=\frac{\phi_S\left(r, t_{n+1}^{-}\right)}{\left\|\phi_S\left(r, t_{n+1}^{-}\right)\right\|_{L^2(\mathbb{R}^+)}}, \quad n \geq 0, \\
&\phi_S\left(r, 0\right) = \phi_S^0(r), \quad \text{with} \quad  \|\phi_S^0\|^2_{L^2(\mathbb{R}^+)} =1,
\end{split}
\end{align}
where $\phi_S\left(r, t_{n}^{\pm}\right) = \lim_{t\rightarrow t_n^{\pm}} \phi_S\left(r, t\right)$, and $\phi_S^0(r)$ is an initial guess for the symmetric or central vortex state solution.

GFDN \eqref{GFDN} can be viewed as the first-order splitting method for the following continuous normalized gradient flow  (CNGF) method \cite{baoComputing}:
\begin{align} \label{CNGF}
\begin{split}
&\p_t \phi_S = \frac12\triangle_{r,S}\phi_{S}+\big(\beta\phi_{S}^2+H\big)\phi_{S}+ \mu_{\phi_S} (t) \phi_S, \quad t_{n}<t<t_{n+1}, \quad n \geq 0, \\
&\phi_S\left(r, 0\right) = \phi_S^0(r), \quad \text{with} \quad  \|\phi_S^0\|^2_{L^2(\mathbb{R}^+)} =1,
\end{split}
\end{align}
where
\begin{align}
\mu_{\phi_S} (t) = \frac{2\pi}{\left\|\phi_S\left(\cdot, t\right)\right\|_{L^2(\mathbb{R}^+)}^2}\int_0^{\infty}\bigg[\frac12\left(\left(\phi_{S}^{\prime}\right)^{2} +\frac{S^2}{r^2} \phi_{S}^{2}\right) -\beta\phi_{S}^4-H \phi_{S}^2\bigg]\ rdr.
\end{align}
It is proved that CNGF \eqref{CNGF} is normalization-conservative and energy-diminishing \cite{baoComputing}.

To improve the GFDN approach to avoiding converging to a "wrong" steady state solution  \cite{MR4214366}  and to using a global strategy for the numerical computation of the local minima,  we devise the  following accelerated, stabilizer-based normalized gradient flow algorithm with Lagrange multipliers (ASGF) to compute the symmetric and central vortex state numerically. We add a nonlocal inertia term and modify the relaxation time in the GFDN model in \eqref{GFDN} as follows
\begin{align} \label{GFLM with stibilizer}
\begin{split}
&\bigg(\alpha_0+\big(\alpha_1-\alpha_2\triangle_{r,S}\big)\p_t\bigg) \p_t \phi_S = \frac12\triangle_{r,S}\phi_{S}+\big(\beta\phi_{S}^2+H\big)\phi_{S}+\mu_{\phi_S}(t_n) \phi_S(r, t_n), \quad t_{n}<t<t_{n+1}, \quad n \geq 0, \\
&\phi_S\left(r, t_{n+1}\right) \triangleq \phi_S\left(r, t_{n+1}^{+}\right)=\frac{\phi_S\left(r, t_{n+1}^{-}\right)}{\left\|\phi_S\left(r, t_{n+1}^{-}\right)\right\|_{L^2(\mathbb{R}^+)}}, \quad n \geq 0,\\
&\phi_S\left(r, 0\right) = \phi_S^0(r), \quad \text{with} \quad  \|\phi_S^0\|^2_{L^2(\mathbb{R}^+)} =1,
\end{split}
\end{align}
where the inertia
\begin{align}
\alpha_0+\big(\alpha_1-\alpha_2\triangle_{r,S}\big)\p_t, \quad \alpha_0, \alpha_1, \alpha_2\geq 0,
\end{align}
serves as a stabilizer, and
\begin{align}
\mu_{\phi_S}(t_n)= \mu \big(\phi_S(\cdot,t_n)\big) = 2\pi\int_0^{\infty}\bigg[\frac12\left(\left(\phi_{S}^{\prime}(r, t_n)\right)^{2} +\frac{S^2}{r^2} \phi_{S}^{2}(r, t_n)\right) -\beta\phi_{S}^4(r, t_n)-H(r, t_n) \phi_{S}^2(r, t_n)\bigg]\ rdr.
\end{align}

\begin{remark}
When $\alpha_0=1, \alpha_1=\alpha_2=0$, this is exactly the GFLM method used in \cite{MR4214366}.
Considering the following continuous Fourier wave in the polar coordinate
\begin{align}
\phi_S(r, \theta, t) = e^{i(k_x\cdot r cos\theta + k_y\cdot r sin\theta + \omega t)}, \quad \vec{k} = (k_x, k_y)^T,
\end{align}
and plugging it into the linear part $\bigg(\alpha_0+\big(\alpha_1-\alpha_2\triangle_{r,\theta}\big)\p_t\bigg) \p_t \phi_S = \frac12\triangle_{r,\theta}\phi_{S}$, where $\triangle_{r,\theta} \triangleq \p_r^2+\frac1{r}\p_r+\frac{1}{r^2} \p_{\theta}^2$, we have
\begin{align}
\omega=\frac{i \alpha_0 \pm \sqrt{-\alpha_0^2+2|\vec{k}|^2 \cdot \big(\alpha_1+\alpha_2 |\vec{k}|^{2}\big)}}{2\big(\alpha_1+\alpha_2 |\vec{k}|^{2}\big)},
\end{align}
from which we can see a smaller value of $\alpha_0$ combined with larger values of $\alpha_1$ and $\alpha_2$ leads to weakened damping of oscillations of the inertia-augmented system.
\end{remark}

Denote $\dot{\phi}_S\triangleq\p_t \phi_S$. Eq. \eqref{GFLM with stibilizer}  can be rewritten as the following gradient flow system
\begin{align}\label{GFLM with stibilizer 1}
\begin{split}
&\p_t \phi_S=\dot{\phi}_S,\\
&\bigg(\alpha_0+\big(\alpha_1-\alpha_2 \triangle_{r,S}\big)\p_t\bigg) \dot{\phi}_S  = \frac12\triangle_{r,S}\phi_{S}+\big(\beta\phi_{S}^2+H\big)\phi_{S}+\mu_{\phi_S}(t_n) \phi_S(r, t_n),\\
&\phi_S\left(r, t_{n+1}\right) \triangleq \phi_S\left(r, t_{n+1}^{+}\right)=\frac{\phi_S\left(r, t_{n+1}^{-}\right)}{\left\|\phi_S\left(r, t_{n+1}^{-}\right)\right\|_{L^2(\mathbb{R}^+)}}, \quad n \geq 0,\\
&\phi_S\left(r, 0\right) = \phi_S^0(r), \quad \text{with} \quad  \|\phi_S^0\|^2_{L^2(\mathbb{R}^+)} =1, \quad \text{and} \quad \dot{\phi}_S\left(r, 0\right) = \dot{\phi}_S^0(r).
\end{split}
\end{align}

\subsubsection{Spatial discretization}

 We map $[0,R]$ into $[-1,1]$ using transformation $r=\frac{R}2(x+1)$, where $x\in I\triangleq[-1,1]$ and denote $u(x)=\phi_S(\frac{R}2(x+1))$, $v(x)=\dot{\phi}_S(\frac{R}2(x+1))$, and $\tilde{H}(x)=H(r)$. We use the Legendre-Galerkin method in $x\in [-1,1]$. Then, \eqref{GFLM with stibilizer 1} is rewritten into the following for $x\in [-1,1]$:
\begin{align}\label{scaling BFET}
\begin{split}
&\p_t u=v, \quad t_{n}<t<t_{n+1}, \quad n \geq 0,\\
&\bigg(\alpha_0+\big(\alpha_1-\frac{4\alpha_2}{R^2} \triangle_{x,S}\big)\p_t\bigg) v  = \frac{2}{R^2}\triangle_{x,S} u +\big(\beta u^2+\tilde{H}\big) u+\mu\big(u(t_n)\big) u(t_n), \quad t_{n}<t<t_{n+1}, \quad n \geq 0,\\
&u(x, t_{n+1})\triangleq\frac{u(x, t^+_{n+1})}{\|u(x, t^+_{n+1})\|_{L^2(I)}}, \quad n \geq 0,\\
& u^0(x) = \phi_S^0(r), \quad v^0(x) =\dot{\phi}_S^0(r),
\end{split}
\end{align}
where  $\|f\|^2_{L^2(I)}=\frac{\pi R^2}2 \int_{-1}^{1}f^2(x)(x+1)dx$, and
\begin{align}
\triangle_{x, S}\triangleq\p_x^2+\frac1{x+1}\p_x-\frac{S^2}{(x+1)^2}.
\end{align}
Now, chemical potential $\mu$ at time $t=t_n$ is rewritten into
\begin{align}
\mu(u(t_n))=&\pi\int_{-1}^{1}\bigg[\big(\p_x u(t_n)\big)^2+\frac{S^2}{(x+1)^2} u^{2}(t_n)-\frac{R^2}{2}\big( \beta u^2(t_n)+\tilde{H}(t_n)\big) u^2(t_n)\bigg]\ (x+1)dx
\end{align}

Given an integer N, we choose $P_N$ from the space of polynomials of degree less than or equal to N, and define
\begin{align}
X_N(S)=\{u\in P_N: u(\pm 1)=0\} \quad \text{for} \quad S\neq 0, \quad X_N(0)=\{u\in P_N: u( 1)=0\}.
\end{align}
Then, we consider the following Legendre-Galerkin approximation to Eq. \eqref{scaling BFET}, where the weight function (x+1) is the Jacobian of the polar transformation. We search for $(u_N, v_N) \in X_N(S)$ such that
\begin{align}\label{scaling BFET integration}
\begin{split}
&\int_{I}\p_t u_N \cdot\omega (x+1)dx=\int_{I}v_N \cdot\omega (x+1)dx,\quad t_{n}<t<t_{n+1}, \quad n \geq 0,\\
&\bigg(\alpha_0+\big(\alpha_1-\frac{4\alpha_2}{R^2} \triangle_{x,S}\big)\p_t\bigg) v_N\cdot\omega (x+1)dx \\
&= \int_{I} \frac{2}{R^2}\triangle_{x,S} u_N\cdot\omega (x+1)dx +\int_{I}I_N g(u)\cdot\omega (x+1)dx+\mu(u_N(t_n)) \int_{I}u_N(t_n)\cdot\omega (x+1)dx, \\
& t_{n}<t<t_{n+1}, \quad n \geq 0,\\
&u_N(t_{n+1})\triangleq\frac{u_N(t^+_{n+1})}{\|u_N(t^+_{n+1})\|_{L^2(I)}}, \quad n \geq 0, \quad \forall \omega \in X_N(S),
\end{split}
\end{align}
where $g(u)=\big(\beta u^2+\tilde{H}\big)u$ and $I_N f$ is the interpolation of $f$ in $\mathbb{P}_N$ at Legendre-Gauss-Lobotta collocation points.

To make the solution satisfying the boundary condition when $S\neq 0$, we construct the following function space:
\begin{align}
X_N(S)=span\{\chi_i(x)=L_i(x)-L_{i+2}(x), \quad i=0,\cdots, N-2\},
\end{align}
where $L_i(x)$ is the $i$ th degree Legendre polynomial.
We define
\begin{align}
\begin{split}
&a_{i j}= \int_{I}\chi_{j}^{\prime} \chi_{i}^{\prime} (x+1) d x, \quad A=\left(a_{i j}\right)_{i, j=0,1, \ldots, N-2}, \\
&b_{i j}= \int_{I} \frac{1}{x+1} \chi_{j} \chi_{i} d x,  \quad B=\left(B_{i j}\right)_{i, j=0,1, \ldots, N-2}, \\
&c_{i j}= \int_{I} \chi_{j} \chi_{i} (x+1) d x, \quad C=\left(C_{i j}\right)_{i, j=0,1, \ldots, N-2}, \\
&u_N=\sum_{i=0}^{N-2} \hat{u}_{i} \chi_{i}(x), \quad\vec{\hat{u}} =(\hat{u}_{0} , \cdots, \hat{u}_{N-2})^T,\\
&v_N=\sum_{i=0}^{N-2} \hat{v}_{i}\chi_{i}(x),\quad \vec{\hat{v}} =(\hat{v}_{0} , \cdots, \hat{v}_{N-2})^T\\
&\big(I_N g(u)\big)(x)= \sum_{i=0}^{N-2} \hat{g}_{i} \chi_{i}(x), \quad\vec{\hat{g}}=(\hat{g}_{0} , \cdots, \hat{g}_{N-2})^T.
\end{split}
\end{align}

 The following results follow from the orthogonality  of the Legendre polynomials.
\begin{lemma}
	Matrix A and B are symmetric, tri-diagonal and given by
	\begin{align}
	a_{i j}=\left\{\begin{array}{ll}2 i+4, & j=i+1, \\ 4 i+6, & j=i,\end{array} \quad b_{i j}=\left\{\begin{array}{ll}-\frac{2}{i+2}, & j=i+1, \\ \frac{2(2 i+3)}{(i+1)(i+2)}, & j=i.\end{array}\right.\right.
	\end{align}
Matrix C is symmetric and seven-diagonal with
	\begin{align}
	c_{i j}=\left\{\begin{array}{ll}
	-\frac{2(i+3)}{(2 i+5)(2 i+7)}, & j=i+3, \\
	-\frac{2}{2 i+5}, & j=i+2, \\
	\frac{2}{(2 i+1)(2 i+5)}+\frac{2(i+3)}{(2 i+5)(2 i+7)}, & j=i+1, \\
	\frac{2}{2 i+1}+\frac{2}{2 i+5}, & j=i.
	\end{array}\right.
	\end{align}
\end{lemma}

Eq. \eqref{scaling BFET integration} with $\omega=\chi_i(x)$ reduces to
\begin{align}\label{eq in frequency}
\begin{split}
&\p_t \vec{\hat{u}}=\vec{\hat{v}},\quad t_{n}<t<t_{n+1}, \quad n \geq 0,\\
&\big[(\alpha_0+\alpha_1\p_t)C + \frac{4\alpha_2}{R^2}(A+S^2 B)\p_t\big]\vec{\hat{v}}= -\frac2{R^2}(A+S^2 B)\vec{\hat{u}}+ C\big(\vec{\hat{g}}+ \mu(u_N^n)\vec{\hat{u}}^n\big),\quad t_{n}<t<t_{n+1}, \quad n \geq 0,\\
&u_N(t_{n+1})\triangleq\frac{u_N(t^+_{n+1})}{\|u_N(t^+_{n+1})\|_{L^2(I)}}, \quad n \geq 0,
\end{split}
\end{align}

In the case when $S = 0$, we construct the following function space:
\begin{align}
X_N(0)=span\{\chi_i(x)=L_i(x)-L_{i+1}(x), \quad i=0,\cdots, N-1\}.
\end{align}

Defining
\begin{align}
\begin{split}
&a_{i j}= \int_{I}\chi_{j}^{\prime} \chi_{i}^{\prime} (x+1) d x, \quad A=\left(a_{i j}\right)_{i, j=0,1, \ldots, N-1}, \\
&c_{i j}= \int_{I}\chi_{j} \chi_{i} (x+1) d x, \quad C=\left(C_{i j}\right)_{i, j=0,1, \ldots, N-1}, \\
&u_N=\sum_{i=0}^{N-1} \hat{u}_{i} \chi_{i}(x),\quad\vec{\hat{u}} =(\hat{u}_{0} , \cdots, \hat{u}_{N-1})^T, \\
&v_N=\sum_{i=0}^{N-1} \hat{v}_{i} \chi_{i}(x),\quad\vec{\hat{v}}=(\hat{v}_{0} , \cdots, \hat{v}_{N-1})^T,\\
&\big(I_N g(u)\big)(x)= \sum_{i=0}^{N-1} \hat{g}_{i} \chi_{i}(x), \quad\vec{\hat{g}}=(\hat{g}_{0} , \cdots, \hat{g}_{N-1})^T,
\end{split}
\end{align}
 one obtains the following results.
\begin{lemma}
	Matrix A is diagonal with
	\begin{align}
	a_{ii} = 2i + 2.
	\end{align}
Matrix C is symmetric and penta-diagonal with
	\begin{align}
	c_{i j}=\left\{\begin{array}{ll}
	-\frac{2(i+2)}{(2 i+3)(2 i+5)}, \quad & j=i+2, \\
	 \frac{4}{(2 i+1)(2 i+3)(2 i+5)}, & j=i+1, \\
	\frac{4(i+1)}{(2 i+1)(2 i+3)}, & j=i.
	\end{array}\right.
	\end{align}
\end{lemma}

Then, eq. \eqref{scaling BFET integration} with $\omega=\chi_i(x)$ reduces to
\begin{align}\label{eq in frequency S=0}
\begin{split}
&\p_t \vec{\hat{u}}=\vec{\hat{v}},\quad t_{n}<t<t_{n+1}, \quad n \geq 0,\\
&\big[(\alpha_0+\alpha_1\p_t)C + \frac{4\alpha_2}{R^2}A\p_t\big]\vec{\hat{v}}= -\frac2{R^2}A\vec{\hat{u}}+ C\big(\vec{\hat{g}}+ \mu(u_N^n)\vec{\hat{u}}^n\big),\quad t_{n}<t<t_{n+1}, \quad n \geq 0,\\
&u_N(t_{n+1})\triangleq\frac{u_N(t^+_{n+1})}{\|u_N(t^+_{n+1})\|_{L^2(I)}}, \quad n \geq 0.
\end{split}
\end{align}

In the following, we  address the issue of solving the transformed magnetic field $\tilde{H}$ equation by applying the Legendre-Galerkin method. In the  transformed coordinate, eq. \eqref{poisson} with the boundary conditions is rewritten into
\begin{align}\label{poisson coordinate transformation}
\begin{split}
-\frac{4}{R^2}\big(\frac{d^2 \tilde{H}}{d x^2} + \frac1{x+1} \frac{d \tilde{H}}{d x}\big) &= \gamma u^2, \quad x \in [-1,1],\\
\frac{d \tilde{H}}{d x}\bigg|_{x=-1}&=0, \quad
\frac{d \tilde{H}}{d x}\bigg|_{x=1}=\frac{\tilde{H}(1)}{2 ln R}.
\end{split}
\end{align}
We seek an approximation of $\tilde{H}$ in space
\begin{align}
Y_N=\left\{\tilde{H} \in \mathbb{P}_{N}: \frac{d\tilde{H}}{dx}\bigg|_{x=-1}=0, \quad  \frac{d\tilde{H}}{dx}\bigg|_{x=1}-\frac{\tilde{H}(1)}{2 ln(R)} =0\right\}.
\end{align}

We define basis functions as follows
\begin{align}
\zeta_i(x)=L_i(x)+a_i L_{i+1}(x) +b_i L_{i+2}(x), \quad i=0, 1, \cdots, N-1,
\end{align}
where $a_i$ and $b_i$ are such that $\zeta_i(x)$ satisfies the  boundary conditions of the function in $Y_N$.

Solving the linear algebra equations, we obtain
\begin{align}
a_i = \frac{2 i+3}{(i+2)^2 \big(ln R (i+1) (i+3) -1\big)}, \quad b_i = \frac{(i+1)^2}{(i+2)^2} \cdot \frac{ln R\cdot i(i+2)-1}{1-ln R (i+1)(i+3)}, \quad i=0, 1, \cdots, N-1.
\end{align}

The Legendre-Galerkin approximation to eq. \eqref{poisson coordinate transformation} is equivalent to  finding $\tilde{H}_N=\sum_{j=0}^{N-2}\hat{\tilde{H}}_j\zeta_j(x) \in Y_N$ such that
\begin{align}
-\frac{4}{R^2}\sum_{j=0}^{N-2}\int_{-1}^1\big(\frac{d^2 \zeta_j}{d x^2} + \frac1{x+1} \frac{d \zeta_j}{d x}\big) \zeta_i(x) (x+1) dx \hat{\tilde{H}}_j= \int_{-1}^1 I_N (\gamma u^2) \zeta_i(x) (x+1) dx.
\end{align}
Let
\begin{align}
h_{i,j}=\int_{-1}^1\big(\frac{d^2 \zeta_j}{d x^2} + \frac1{x+1} \frac{d \zeta_j}{d x}\big) \ \zeta_i(x) (x+1) dx,\quad H_P=\left(h_{i j}\right)_{i, j=0,1, \ldots, N-1}.
\end{align}
We have the following lemma.
\begin{lemma}
	Matrix $H_P$ is symmetric and tri-diagonal with
	\begin{align}
	h_{i j}=\left\{\begin{array}{ll}2(i+4) a_i, \quad & j=i+1 \\
	 2(i+1) a_i +2 (2 i+3) b_i + 2(i+2) a_i b_i, & j=i\end{array}\right.
	\end{align}
\end{lemma}

In the spectral representation of solutions, we define
\begin{align}
\begin{split}
&u_N^n=\sum_{i=0}^{N_1} \hat{u}_{i}^n \chi_{i}(x),\quad\vec{\hat{u}}^{n} =(\hat{u}_{0}^n , \cdots, \hat{u}_{N_1}^n)^T, \quad u_N^*=\sum_{i=0}^{N_1} \hat{u}_{i}^* \chi_{i}(x),\quad\vec{\hat{u}}^{*}=(\hat{u}_{0}^* , \cdots, \hat{u}_{N_1}^*)^T,\\
&v_N^n=\sum_{i=0}^{N_1} \hat{v}_{i}^n \chi_{i}(x),\quad\vec{\hat{v}}^{n} =(\hat{v}_{0}^n , \cdots, \hat{v}_{N_1}^n)^T,\quad v_N^*=\sum_{i=0}^{N_1} \hat{v}_{i}^* \chi_{i}(x),\quad\vec{\hat{v}}^{*}=(\hat{v}_{0}^* , \cdots, \hat{v}_{N_1}^*)^T,\\
&\big(I_N g(u^n)\big)(x)= \sum_{i=0}^{N1} \hat{g}_{i}^n \chi_{i}(x), \quad\vec{\hat{g}}^n=(\hat{g}_{0}^n , \cdots, \hat{g}_{N1}^n)^T, \\
&\text{with} \ N1=N-2\ \text{when} \ S\neq 0, \quad N1=N-1 \ \text{when} \  S= 0.
\end{split}
\end{align}

\subsubsection{Temporal discretization}
After the Legendre-Galerkin approximation in space, we present two time-discredized schemes below. 
We use backward and forward mixed Euler scheme to discretize ODE system   \eqref{eq in frequency} for $S>0$ or \eqref{eq in frequency S=0} for $S=0$ in time to arrive at the first ASGF algorithm.
\begin{algorithm}[\bf $ASGF-\uppercase\expandafter{\romannumeral 1\relax}$]

Given initial data $u^0, v^0$, compute the spectral coefficients $\vec {\hat{u}}^0, \vec{\hat{v}}^0$. For $n>0$, compute $\vec{\hat{u}}^n, \vec{\hat{v}}^n$ via
\begin{align}\label{ASGF 1}
\begin{split}
&\frac{\vec{\hat{u}}^{*}-\vec{\hat{u}}^{n}}{\triangle t}=\vec{\hat{v}}^{*} ,\quad t_{n}<t<t_{n+1}, \quad n \geq 0, \\
&\alpha_0C\vec{\hat{v}}^{*}  + \big(\alpha_1C+\frac{4\alpha_2}{R^2}(A+S^2 B)\big) \frac{\vec{\hat{v}}^{*}-\vec{\hat{v}}^{n}}{\triangle t} = \big(-\frac2{R^2}(A+S^2 B)-\alpha C\big)\vec{\hat{u}}^{*}+ C\big(\alpha\vec{\hat{u}}^{n} +\vec{\hat{g}}^n+ \mu(u_N^n)\vec{\hat{u}}^{n}\big),\\
&\quad t_{n}<t<t_{n+1}, \quad n \geq 0,\\
&u_N^{n+1}=\frac{u_N^{*}}{\|u_N^{*}\|_{L^2(I)}},\quad n \geq 0, \quad \forall S \geq 0,
\end{split}
\end{align}
where $\alpha = max\{-\frac12\big(\beta |u^n|^2+ \tilde{H}^n +\mu(u^n)\big),0\}$ is a chosen stabilization parameter such that the time step can be as large as possible;

until the $L_2$-norm of the residue is less than tolerance $\epsilon$.
\end{algorithm}

In order to preserve the gradient-flow structure at the discrete level,  different u-dependent terms on the right-hand side in the second expression of \eqref{eq in frequency}  and \eqref{eq in frequency S=0} should be approximated at the same time level, i.e., explicit  treatment.
Hence, we use the combined backward Euler method on the first equation of \eqref{eq in frequency} or \eqref{eq in frequency S=0} and the forward Euler method on the second equation, we end up with the second scheme as follows.
\begin{algorithm}[\bf $ASGF-\uppercase\expandafter{\romannumeral 2\relax}$] Given initial data $u^0, v^0$, compute the spectral coefficients $\vec {\hat{u}}^0, \vec{\hat{v}}^0$. For $n>0$, compute $\vec{\hat{u}}^n, \vec{\hat{v}}^n$ via
\begin{align}\label{ASGF 2}
\begin{split}
&\frac{\vec{\hat{u}}^{*}-\vec{\hat{u}}^{n}}{\triangle t}=\vec{\hat{v}}^{*} ,\quad t_{n}<t<t_{n+1}, \quad n \geq 0, \\
&\alpha_0C\vec{\hat{v}}^{*}  + \big(\alpha_1C+\frac{4\alpha_2}{R^2}(A+S^2 B)\big) \frac{\vec{\hat{v}}^{*}-\vec{\hat{v}}^{n}}{\triangle t} = -\frac2{R^2}(A+S^2 B)\vec{\hat{u}}^{n}+ C\big(\vec{\hat{g}}^n+ \mu(u_N^n)\vec{\hat{u}}^{n}\big),\\
&\quad t_{n}<t<t_{n+1}, \quad n \geq 0,\\
&u_N^{n+1}=\frac{u_N^{*}}{\|u_N^{*}\|_{L^2(I)}},\quad n \geq 0, \quad \forall S \geq 0;
\end{split}
\end{align}
until the $L_2$-norm of the residue is less than tolerance $\epsilon$.
\end{algorithm}

Note that both schemes are first order in time. Next, we compare the new  numerical schemes $ASGF-\uppercase\expandafter{\romannumeral 1\relax}$, $ASGF-\uppercase\expandafter{\romannumeral 2\relax}$ with the existing GFLM method, which is the special case of the above two schemes with $\alpha_0=1, \alpha_1=\alpha_2=0$, when solving central vortex state solutions.

\subsubsection{Numerical results of the single-component GP model}

\begin{figure}
	\begin{minipage}[t]{0.5\linewidth}
		\centering     
		\includegraphics[width=\textwidth]{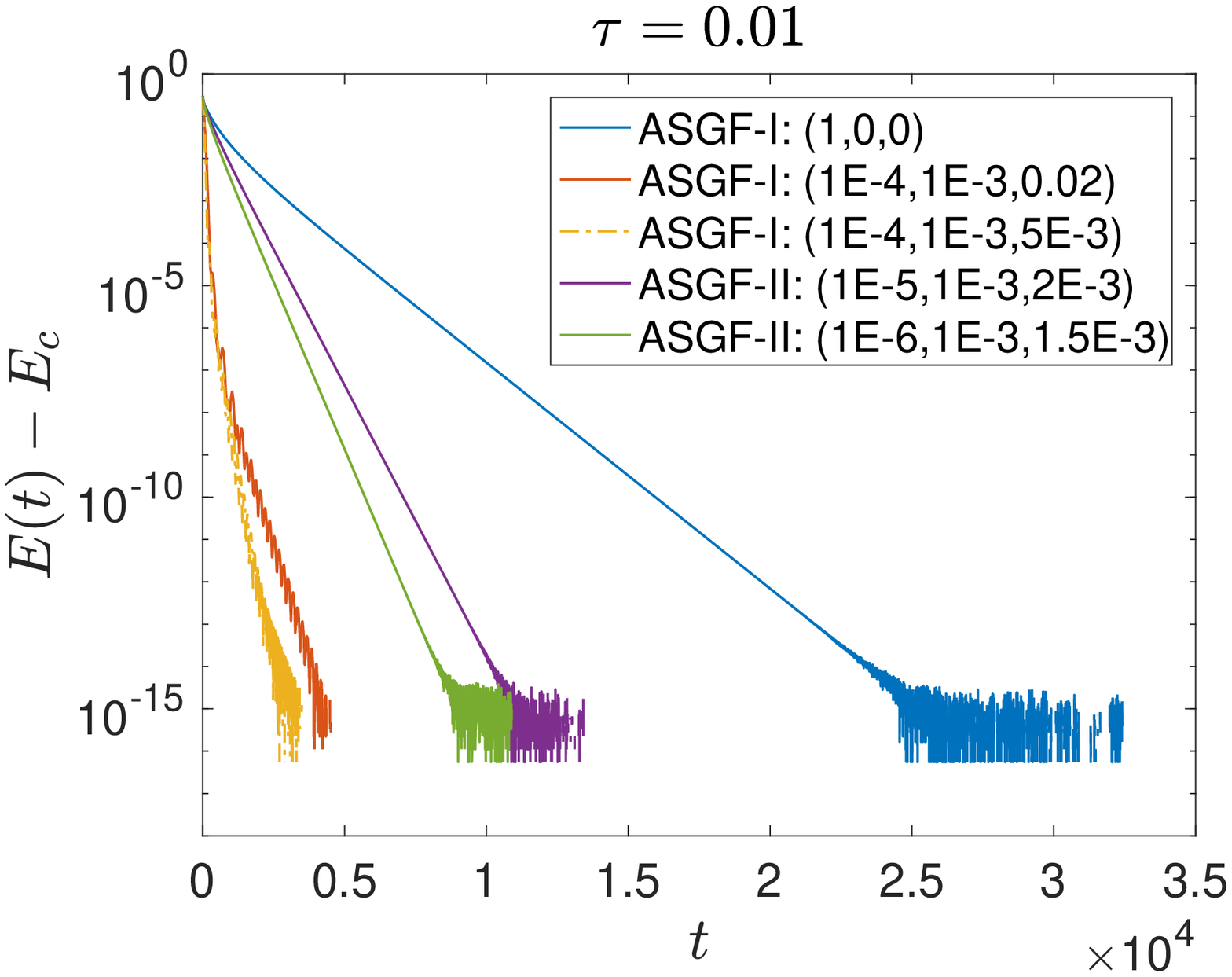}
	\end{minipage} 
	\hfill
	\begin{minipage}[t]{0.5\linewidth}
		\centering
		\includegraphics[width=\textwidth]{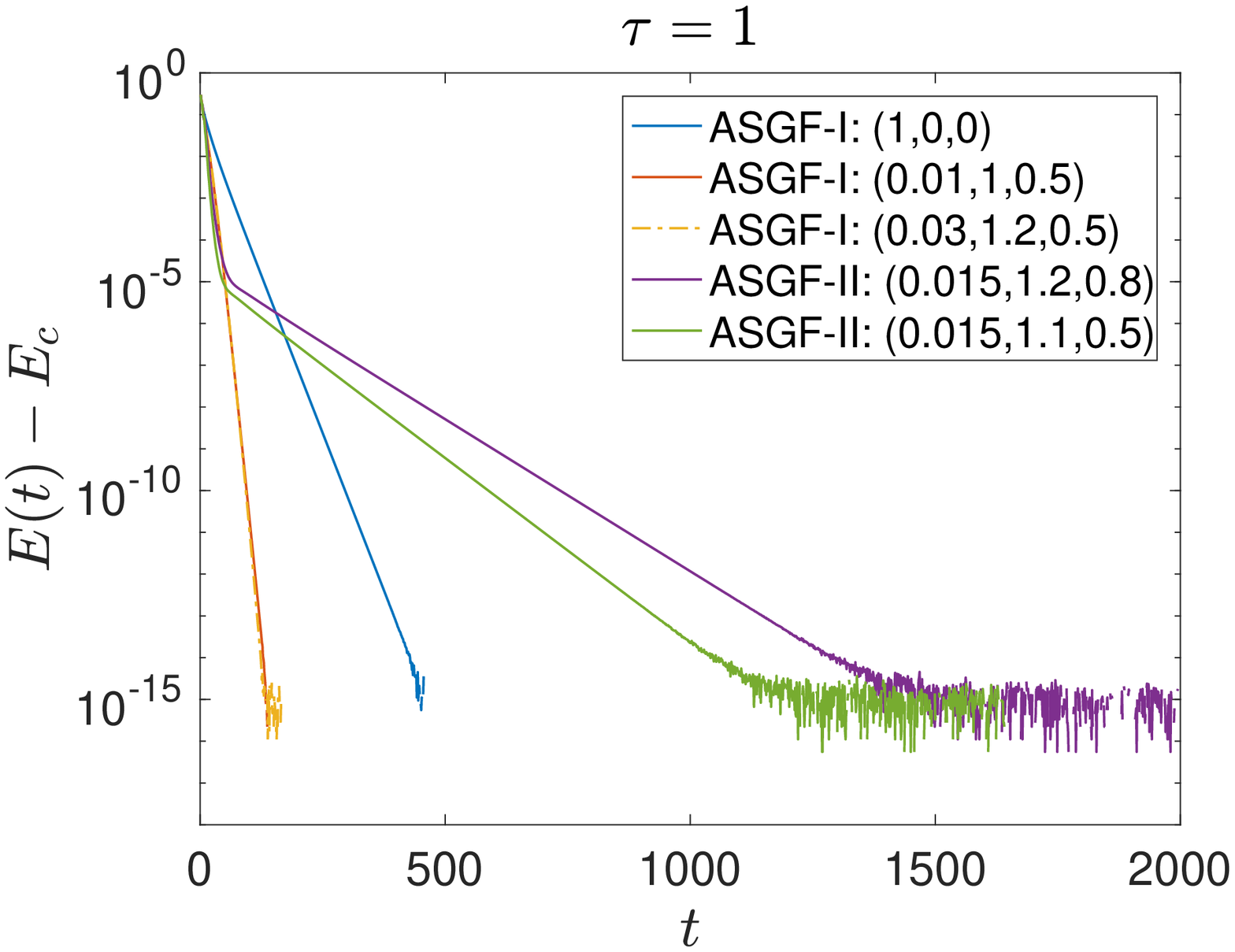}
	\end{minipage}
	\caption{Time evolution of relative energy $E(\phi_S(\cdot, t)) - E_c$ (in logarithmic scale) by different numerical schemes when computing the central vortex state solution. Left: $\tau =0.01$; Right: $\tau =1$. The triplet, for example $(1,0,0)$, in the inset represents $\alpha_{0}=1, \alpha_{1}=0, \alpha_{2}=0$, respectively.}
	\label{Fig21}
\end{figure}

The initial condition for the iterative scheme is chosen as $\phi_S^0(r) = \frac 1{\sqrt{\pi S!}} r^S e^{-r^2/2}$ and $\dot{\phi}_S^0 = 0$. The stopping criterion for the time marching is that the $L^2$-norm of residue of the Euler–Lagrange equation \eqref{single chemical polar} is less than  given tolerance $\varepsilon = 10^{-10}$.
We take $S=2$, $\beta=30$, $\gamma=\pi$, $U=[0,20]$ with the dimension of the discrete Legendre space $N_l = 200$. The performance of two methods with  a few selected time step $\tau$ are tabulated  in Table \ref{Tab21}, where $E_c$ and $\mu_c$ are the energy and chemical potential at the central vortex state solution, respectively, and \#$iter$ is the number of iterations (or time marching steps) used in computing the steady state. The time evolution of relative energy $E(\phi_S(\cdot, t)) - E_c$ (in logarithmic scale) by the use of  the two numerical schemes for computing the central vortex state solutions are shown in Figure \ref{Fig21}.

From the numerical results obtained using the $ASGF-\uppercase\expandafter{\romannumeral 1\relax}$ scheme in Tab. \ref{Tab21} and Fig. \ref{Fig21}, we observe that  the computational time of the GFLM method, corresponding to $\alpha_0=1, \alpha_1=\alpha_2=0$, is much longer (about six times) than that of the ASGF method for the same time step, demonstrating that the $ASGF-\uppercase\expandafter{\romannumeral 1\relax}$ scheme is more efficient. Hence,  $\alpha_0+\big(\alpha_1-\alpha_2 \triangle_{r,S}\big)\p_t$ indeed speeds up convergence to the steady state of the normalized gradient flow solution. Likewise, the $ASGF-\uppercase\expandafter{\romannumeral 2\relax}$ scheme outperforms the GFLM scheme considerably as well thanks to stabilizer $\alpha_0+\big(\alpha_1-\alpha_2 \triangle_{r,S}\big)\p_t$ again. We thus conclude that both ASGF schemes have better convergence properties than the GFLM scheme.

\begin{table}
	\centering
	\begin{tabular}{|c|c|ccccccc|}
		\hline
		Scheme&$\tau$ &$\alpha_0$& $\alpha_1$&$\alpha_2$&$\mathrm{CPU}(\mathrm{s})$ & $ E_{c}$ & $\mu_{c} $& \#$iter$ \\
		\hline
		\multirow{9}{*}{ASGF-\uppercase\expandafter{\romannumeral 1\relax}}& &1&0&0&323.5&0.4666956706&0.5688732593&32431 \\
		&0.01&0.0001&0.001&0.01&74.8&0.4666956706&0.5688732599&7199 \\
		&&0.0001&0.001&0.005&50.1&0.4666956706&0.5688732599 &4785\\
		\cline{2-9}
		&&1&0&0&37.5&0.4666956706&0.5688732593&3487\\
		&0.1&0.001&0.01&0.05&7.8&0.4666956706&0.5688732599&711\\
		&&0.01&0.01&0.05&5.2&0.4666956706&0.5688732599&475\\
		\cline{2-9}
		&&1&0&0&6.60&0.4666956706&0.5688732594&590 \\
		&1&0.01&1&0.5&2.3&0.4666956706&0.5688732600&213 \\
		&&0.03&1.2&0.5&1.9&0.4666956706&0.5688732600&168 \\
		\hline
		\multirow{9}{*}{ASGF-\uppercase\expandafter{\romannumeral 2\relax}}
		&&1&0&0&\quad-&\quad-&\quad -&\quad-\\
		&0.01&1E-5&0.001&0.002&139.1&0.4666956706&0.5688732599&13448 \\
		&&1E-6&0.001&0.0015&118.9&0.4666956706&0.5688732599&10889\\
		\cline{2-9}
		&&1&0&0&\quad-&\quad \quad -&\quad \quad -&\quad - \\
		&0.1&0.001&0.01&0.05&88.9&0.4666956706&0.5688732598&8376\\
		&&0.0015&0.01&0.02&30.3&0.4666956706&0.5688732596&2873\\
		\cline{2-9}
		&&1&0&0&\quad-&\quad \quad -&\quad \quad -&\quad - \\
		&1&0.015&1.2&0.8&21.2&0.4666956706&0.5688732596&2007 \\
		&&0.015&1.1&0.5&18.5&0.4666956706&0.5688732600&1642 \\
		\hline
	\end{tabular}
	\caption{Performance of the  ASGF-\uppercase\expandafter{\romannumeral 1\relax} scheme and the \uppercase\expandafter{\romannumeral 2\relax} scheme when computing the  central vortex state solution  of a nonlocal single component GP model. The  ASGF-\uppercase\expandafter{\romannumeral 1\relax} scheme outperforms the \uppercase\expandafter{\romannumeral 2\relax} scheme, and both are better than the GFLM scheme. }
	\label{Tab21}
\end{table}

\begin{remark}
	We have conducted additional numerical simulations for the rotating GP equation in a separate study, such as the logarithmic Schr\"odinger equation with the angular momentum, and  noted the superior performance of the ASGF approach than that of GFLM. In fact, the GFLM scheme can only converge if time step $\triangle t$ is less than or equal to $\mathcal{O}(10^{-3})$ while  the steady states can be reached with the ASGF scheme even at $\triangle t \thicksim \mathcal{O}(1)$.
\end{remark}

\subsection{Perturbed preconditioned nonlinear conjugate gradient (PPNCG) method}

Here, we  present a projection  method  which preserves the gradient-flow structure of \eqref{single energy functional_polar} at the discrete level while implicitly accounting for the presence of the unit-norm constraint \eqref{single polar constraint}.
The projected, preconditioned conjugate gradient method (PCG) for the minimization of $E(\phi_{S})$ on $\mathcal{S}_1$ is built on the update given by
\begin{align}\label{alpha_n cg}
\begin{split}
\tilde{\phi_{S}}^{n+1} = \phi_{S}^{n}  -\alpha_n P^n r^n,\quad
\phi_S^{n+1}=\frac{\tilde{\phi_{S}}^{n+1} }{\|\tilde{\phi_{S}}^{n+1} \|_{L^2(\mathbb{R}^+)}},
\end{split}
\end{align}
where
\begin{align}
P^n:=\bigg(\pi\int_0^{\infty}\left(\left((\phi_{S}^n)^{\prime}\right)^{2} +\frac{S^2}{r^2} (\phi_{S}^n)^{2}\right)\ rdr - \frac12\triangle_{r,S}\bigg)^{-1}
\end{align}
is the symmetric positive definite preconditioner, and
\begin{align}\label{residual}
 r^n:=-\frac12\triangle_{r,S}\phi_{S}^{n} -\big(\beta (\phi_S^n)^2+H^n\big)\phi_{S}^{n}
  - \mu_n \phi_{S}^{n}
\end{align}
  is the residue.

We reformulate this update formula as follows
\begin{align}\label{theta_n cg}
\phi_S^{n+1} = cos(\theta_n) \phi_S^{n} + sin(\theta_n) \frac{p^n}{\|p^n\|_{L^2(\mathbb{R}^+)}}, \quad p^n = d^n-\langle d^n, \phi_S^{n} \rangle \phi_S^{n},
\end{align}
 where $d^n=-P^n r^n +\beta_n d^{n-1}$ is the search direction, $\langle \cdot, \cdot \rangle $ is the $L^2(\mathbb{R}^+)$ inner product and $\beta_n$ is the "momentum" term chosen to enforce the conjugacy of search directions $d_k, \ k=1, \cdots, n$.   Either one of the following expressions
\begin{align}
\begin{split}
\beta_n=\beta_n^{FR} := \frac{\langle r^n, P^n r^n \rangle}{\langle r^{n-1}, P^{n-1}r^{n-1} \rangle}, \quad \text{(Fletcher–Reeves)},\\
\beta_n=\beta_n^{FR} := max\left\{\frac{\langle r^n-r^{n-1}, P^n r^n \rangle}{\langle r^{n-1}, P^{n-1}r^{n-1} \rangle}, 0\right\}, \quad \text{(Polak–Ribi\'ere))}
\end{split}
\end{align}
can be used to update $\beta_n$.
We note that equation \eqref{alpha_n cg} and \eqref{theta_n cg} are equivalent when $\theta_n$ or $\alpha_n$ is small enough, with a one-to-one correspondence between $\theta_n$ and $\alpha_n$. In practice, $\theta_n$ may not be small, then a general line-minimization approach such as Brent’s algorithm should be adapted.

Expanding $\phi_S^{n+1}$ up to second-order in $\theta_n$, we obtain
\begin{align}
\phi_S^{n+1}=\left(1-\frac{\theta_{n}^{2}}{2}\right) \phi_S^{n}+\theta_{n} \frac{p^{n}}{\left\|p^{n}\right\|}+\mathcal{O}\left(\theta_{n}^{3}\right)
\end{align}
and therefore
\begin{align}
E\left(\phi_S^{n+1}\right)=E\left(\phi_S^{n}\right)+\frac{\theta_{n}}{\left\|p^{n}\right\|} \left\langle\nabla E\left(\phi_S^{n}\right), p^{n}\right\rangle+\frac{1}{2} \frac{\theta_{n}^{2}}{\left\|p^{n}\right\|^{2}}\left[\nabla^{2} E\left(\phi_S^{n}\right)\left[p^{n}, p^{n}\right]-2\mu(\phi_S^n)\left\|p^{n}\right\|^{2}\right]+\mathcal{O}\left(\theta_{n}^{3}\right).
\end{align}
Minimizing the above functional with respect to $\theta_n$  yields
\begin{align}
\theta_{n}^{\mathrm{opt}}=\frac{-\left\langle\nabla E\left(\phi^{n}\right), p^{n}\right\rangle\left\|p^{n}\right\|}{\nabla^{2} E\left(\phi^{n}\right)\left[p^{n}, p^{n}\right]-2\mu(\phi_S^n)\left\|p_{n}\right\|^{2}}.
\end{align}

It's known that the  gradient descent method can be exponentially slow in the presence of saddle points \cite{Du2017}. Due to non-convexity of energy functional \eqref{single energy functional_polar},  its critical points may be an approximate saddle point instead of  a local minimum. Then an appropriate procedure should be put in place to escape from the saddle point. We use the following perturbed preconditioned nonlinear conjugate gradient method (PPNCG) to find the minimum of the Hamiltonian confined in manifold $\mathcal{S}_1$ \cite{ANTOINE201792, Sun2019}.

\begin{algorithm} [PPNCG ]
	\begin{algorithmic}

		\STATE   Initiate $\phi_S^0 \in \mathcal{S}_1$ and use the steepest descent method in the first step.
		\STATE $\mu_0 = \mu(\phi_{S}^{0})$,
		\STATE $r^0=-\frac12\triangle_{r,S}\phi_{S}^{0} -\big(\beta (\phi_S^0)^2+H^0\big)\phi_{S}^{0}
		- \mu_0 \phi_{S}^{0}$,
		\STATE $d^0=-P^0 r^0$, $p^0 = d^0-\langle d^0, \phi_S^{0} \rangle \phi_S^{0}$,
		\STATE $\theta_0 = argmin_{\theta}E\big( cos(\theta) \phi_S^{0} + sin(\theta) \frac{p^0}{\|p^0\|_{L^2(\mathbb{R}^+)}}\big)$,
		\STATE $\phi_S^{1} = cos(\theta_0) \phi_S^{0} + sin(\theta_0)\frac{p^0}{\|p^0\|_{L^2(\mathbb{R}^+)}}$,
		\STATE Set $n=1, \|r^1\| = \|r^0\|$;
 		\WHILE {$\|r^n\| \geq \varepsilon$}
		\STATE  Use the projected conjugate gradient method in this loop.
		\STATE $\mu_n = \mu(\phi_{S}^{n})$,
		\STATE $r^n=-\frac12\triangle_{r,S}\phi_{S}^{n} -\big(\beta (\phi_S^n)^2+H\big)\phi_{S}^{n}
		- \mu_n \phi_{S}^{n}$,
		\STATE $\beta_n= max\left\{\frac{\langle r^n-r^{n-1}, P^n r^n \rangle}{\langle r^{n-1}, P^{n-1}r^{n-1} \rangle}, 0\right\}$,
		\STATE $d^n=-P^n r^n +\beta_n d^{n-1}$,
		\STATE $p^n = d^n-\langle d^n, \phi_S^{n} \rangle \phi_S^{n}$,
		\STATE $\theta_n = argmin_{\theta}E\big(cos(\theta) \phi_S^{n} + sin(\theta) \frac{p^n}{\|p^n\|_{L^2(\mathbb{R}^+)}}\big)$,
		\STATE $\phi_S^{n+1} = cos(\theta_n) \phi_S^{n} + sin(\theta_n) \frac{p^n}{\|p^n\|_{L^2(\mathbb{R}^+)}}$,
		\STATE $n = n+1$;
		\ENDWHILE
		
		\STATE Denote  the above numerical solution by  $\phi_S^{\star}$,		perturb $\phi_S^{\star}$ by adding an appropriate level of noise in its tangent space and  map it back to the manifold, denoted as $\phi_S^{\star,0}$,
		put the perturbed numerical solution $\phi_S^{\star,0}$ into above loop and run a few (about 7) iterations.
		
		\IF {the value of the energy functional decreases}
		 \STATE it indicates that the numerical solution escapes from the approximate saddle point;
		 \ELSE
		  \IF {the value does not decrease}
		 \STATE it is accepted as an approximate  minimum.
		 \ENDIF
		 \ENDIF
		
	\end{algorithmic}

\label{alg: single PPNCG}
\end{algorithm}

\subsubsection{Numerical results of the single-component GP model}

 We compare the computational time and the number of iterations with a few selected values of $S$ and $\beta$ between the ASGF-\uppercase\expandafter{\romannumeral 1\relax} and PPNCG method when computing the symmetric state and central vortex state solution.

The initial condition is chosen as $\phi_S^0(r) = \frac 1{\sqrt{\pi S!}} r^S e^{-r^2/2}$ and the initial velocity $\dot{\phi}_S^0 = 10 \times \phi_S^0(r)$, time step $\tau=1$, $\alpha_0=0.01, \alpha_1=1, \alpha_2=0.2$  for ASGF-\uppercase\expandafter{\romannumeral 1\relax}. The stopping criterion for time marching is that the $L^2$-norm of residue of the Euler–Lagrange equation \eqref{single chemical polar} is less than the given tolerance $\varepsilon = 10^{-10}$.
We take $R=18$ for $S=0$, $R=20$ for $S=2$, $R=30$ for $S=5$ and $R=35$ for $S=8$  with the dimension of the approximate solution space $N_l = 10 R$. The performance of the  two methods are shown in Table \ref{Tab22}, from which we see clearly that the PPNCG method is much better than the ASGF-\uppercase\expandafter{\romannumeral 1\relax} method.

\begin{table}
	\centering
	\begin{tabular}{|cc|ccc|ccc|ccc|ccc|}
		\hline
		\multicolumn{2}{|c|}{S}&&0&& &2&&& 5&&&8& \\
		\hline
		\multicolumn{2}{|c|}{$\beta$} &0&3&4.5&0&30&40&0&50&80&0&100&140\\
		\hline
		\multirow{2}{*}{ASGF-\uppercase\expandafter{\romannumeral 1\relax}} &CPU(s)&1.61&1.48&2.12&2.12&2.48&10.25&12.33&12.99&64.32&24.51&110.72&447.98\\
		&\#$iter$&205&191&272&201&236&985&245&258&1303&274&1148&5104\\
		\hline
		\multirow{2}{*}{PPNCG} &CPU(s)&0.36&0.40&0.39&0.56&0.65&0.68&2.70&3.85&6.06&6.55&13.59&22.44\\
		&\#$iter$&28&32&31&30&39&40&40&56&77&46&95&158\\
		\hline
	\end{tabular}
\caption{ Performance comparison between the ASGF-\uppercase\expandafter{\romannumeral 1\relax} method and the PPNCG method in the total computational time and the number of iterations with respect to various values of $S$ and $\beta$.}
\label{Tab22}
\end{table}

\section{Numerical methods for  the coupled binary GP model}

In this section, we extend the methods developed in \S 2 to  the coupled binary GP model when computing the symmetric state and central vortex state solution.
The  coupled Gross-Pitaevskii equation with  wave function $\Psi=(\psi_{\uparrow},\psi_{\downarrow})^T$ in a dimensionless 2D domain is given by
\begin{align}\begin{split}\label{mGPE}
i\frac{\p \psi_{\uparrow}}{\p t}&=\big[-\frac12\nabla^2+V_{\uparrow}(\xx)-\eta-\beta|\psi_{\downarrow}|^2\big]\psi_{\uparrow}- \breve{H} \psi_{\downarrow},
\\
i\frac{\p \psi_{\downarrow}}{\p t}&=\big[-\frac12\nabla^2+V_{\downarrow}(\xx)+\eta-\beta|\psi_{\uparrow}|^2\big]\psi_{\downarrow} - \overline{\breve{H} } \psi_{\uparrow},
\end{split}
\end{align}
where $\overline{f}$ represents the conjugate of $f$,
$\eta$ is the dimensionless detuning parameter,  $\beta$ is the dimensionless contact interaction parameter, $V_j(\xx) \ (j=1,2)$ are the external potentials, the magnetic field $\breve{H} $ satisfies
\begin{align}\label{H}
-\triangle \breve{H} =\gamma \psi_{\uparrow}(\xx) \overline{\psi_{\downarrow}(\xx)}.
\end{align}
If we use the fundamental solution of the 2D Poisson equation, $\breve{H} $ can be expressed explicitly by
\begin{align}\label{H fundamental solution}
\breve{H}  = H_0 +\breve{H_1} = H_0 -\frac{\gamma }{2 \pi}  \int_{\mathbb{R}^2} \ln \left(\left|\xx-\xx^{\prime}\right|\right)
\psi_{\uparrow}\left(\xx^{\prime}\right) \overline{\psi_{\downarrow}\left(\xx^{\prime}\right)} d \xx^{\prime},
\end{align}
where $H_0$ is a background magnetic field {and $\Delta H_0=0$}.

The Hamiltonian of the conservative system  is identified as
\begin{align}\label{energy functional}
E(\Psi)=&\int_{\mathbb{R}^2}\bigg[\sum_{j=\uparrow,\downarrow}\big(\frac12|\nabla\psi_j|^2+V_{j} |\psi_{j}|^2\big)-\eta\big(|\psi_{\uparrow}|^2-|\psi_{\downarrow}|^2\big)-\beta|\psi_{\uparrow}|^2|\psi_{\downarrow}|^2 \nonumber\\
&- 2H_0  Re{(\psi_{\uparrow} \overline{\psi_{\downarrow}})}-\frac1{\gamma} |\nabla \breve{H_1}|^2 \bigg]d\xx. 
\end{align}

In the following, we assume external potentials $V_{1} = V_{\uparrow}$ and $V_{2}=V_{\downarrow}$ are radially symmetric, and we limit our search for the steady  state of \eqref{energy functional} to the central vortex form in  polar coordinates $(r,\theta)$ as follows
\begin{align}\label{stationary polar}
\psi_{\uparrow}(\xx)=  e^{-i S \theta} \phi_{1}(r), \quad \psi_{\downarrow}(\xx)=  e^{-i S \theta} \phi_{2}(r),
\end{align}
where S is the winding number and $\Phi_{S}:=(\phi_{1}, \phi_{2})^T$ is a real-valued radial wave function vector. The radial steady magnetic field $\breve{H}_1$ can then be expressed as:
\begin{align}\label{binary poisson R+}
-\frac1{r} \frac{d}{dr}\big(r\frac{d H_1}{dr}\big)=\gamma \phi_1 \phi_2, \quad r\in [0,\infty),
\end{align}

The corresponding energy functional \eqref{energy functional} with the radially symmetric solution  reduces to the following, parameterized by winding number S:
\begin{align}\label{binary energy functional_polar}
E(\Phi_S)= &2\pi\int_0^{\infty}\big[\sum_{j=1, 2}\frac12\left(\left(\phi_{j}^{\prime}\right)^{2} +\frac{S^2}{r^2} \phi_{j}^{2}\right) +\sum_{j=1,2}V_j \phi_{j}^2 \nonumber\\
& - \eta\big(\phi_{1}^2-\phi_{2}^2\big)-\beta\phi_{1}^2\phi_{2}^2-(2H_0+H_1) \phi_{1}\phi_{2}\big] rdr,
\end{align}
subject to  constraint
\begin{align}\label{polar constraint}
2\pi\int_0^{\infty}\big(\phi_1^2(r)+\phi_2^2(r)\big)\ r dr=1.
\end{align}
For a given winding number $S$, we denote the symmetric state as $\Phi_S^s$ at $S=0$ and central vortex state as $\Phi_S^c$ at $S>0$, respectively, which minimizes $E(\Phi_S)$ in manifold
\begin{align} \label{binary manifold}
\mathcal{S}_2 :=\left\{\Phi_S=\left(\phi_{1}, \phi_{2}\right)^{T} \ \big|\ \|\Phi_S\|_{L^2(\mathbb{R}^+)}^{2}=2\pi\int_0^{\infty}\big(\phi_1^2(r)+\phi_2^2(r)\big)\ r dr=1, E(\Phi_S)<\infty\right\}.
\end{align}

One deduces the corresponding Euler-Lagrange equations of the constrained minimization problem as follows
\begin{align}\label{binary lagrange equation}
\begin{split}
\mu\phi_{1}&=-\frac12\triangle_{r,S}\phi_{1}+\bigg[V_{1}-\eta-\beta\phi_{2}^2\bigg]\phi_{1}-\big(H_0+H_1\big)\phi_{2},\\
\mu\phi_{2}&=-\frac12\triangle_{r,S}\phi_{2}+\bigg[V_{2}+\eta-\beta\phi_{1}^2\bigg]\phi_{2}-\big(H_0+H_1\big)\phi_{1}.
\end{split}
\end{align}
The corresponding Lagrange multiplier or eigenvalue (chemical potential)  is given by
\begin{align}\label{chemical polar}
\mu(\Phi_S)=&2\pi\int_0^{\infty}\bigg[\sum_{j=1, 2}\frac12\left(\left(\phi_{j}^{\prime}\right)^{2} +\frac{S^2}{r^2} \phi_{j}^{2}\right) +\sum_{j=1,2}V_j \phi_{j}^2-\eta\big(\phi_{1}^2-\phi_{2}^2\big)-2\beta\phi_{1}^2\phi_{2}^2 -2(H_0+H_1) \phi_{1}\phi_{2} \bigg]\ rdr.
\end{align}
We search for steady state  solutions in the function space defined by
\begin{align}
X_S^2(\mathbb{R^+}) :=H^1_S(\mathbb{R^+}) \cap L^4_S(\mathbb{R^+}) \cap L^2_V(\mathbb{R}^+)\cap L_{log,S}^2(\mathbb{R}^+),
\end{align}
where
$L^2_V(\mathbb{R}^+) := \left\{\Phi_S \bigg| \int_0^{\infty} (V_1(r) \phi_1^2(r) +V_2(r) \phi_2^2(r)) r dr < \infty \right\}$, and  $\lim\limits_{r\rightarrow +\infty} V_j(r) = +\infty ( j =1 ,2)$.
\begin{lemma} (see \cite{WANG2021132852})
	In Hilbert space $X_S^2(\mathbb{R}^+)$ and for any given $S\geq 0$,  there exists a symmetric or central vortex steady state $\Phi_S=(\phi_1,\phi_2)^T$ of  \eqref{binary energy functional_polar} when $\beta \leq \beta_b$ and $V_j (r) \geq \frac{\gamma}{8 \pi} r^2, j=1,2$, where $\beta_b$ is defined in Table \ref{beta_b S}.  When $\beta > 2\beta_b$, there does not exist any symmetric or central vortex steady state.
		\begin{table}
		\centering
		\begin{tabular}{p{0.5cm}p{1cm}p{1cm}p{1cm}p{1cm}p{1cm}p{1cm}p{1cm}p{1cm}{c}}
			\hline
			S	& 0 & 1 & 2& 3  & 4 & 5&6&7 \\
			\hline
			$\beta_b$ & 11.70& 48.31&  89.75&132.42&175.50&218.76&262.11&305.51\\
			\hline
			S &8 & 9  & 10 & 11&12 &13 &14 &15 \\
			\hline
			$\beta_b$  &348.94 &392.40&435.87&479.35&522.84&566.33& 609.83&653.34\\
			\hline
		\end{tabular}			
		\caption{$\beta_b$  vs S}
		\label{beta_b S}
	\end{table}
\end{lemma}

\subsection{ASGF method}

The normalized gradient flow model for computing the symmetric and central vortex state of the nonlocal binary GP model reads as follows.
Given the time sequence $t_n = n\tau$ for $n = 0,1,2, \cdots$, for $t \in [t_n, t_{n+1}) (n \geq 0)$, one solves the following gradient flow model in time:
\begin{align} \label{binary GFLM with stibilizer}
\begin{split}
&\bigg(\alpha_{01}+\big(\alpha_{11}-\alpha_{21}\triangle_{r,S}\big)\p_t\bigg) \p_t \phi_1 = \frac12\triangle_{r,S}\phi_{1}-g_1+\mu(\Phi_S(t_n)) \phi_1(r, t_n), \quad t_{n}<t<t_{n+1}, \quad n \geq 0, \\
&\bigg(\alpha_{02}+\big(\alpha_{12}-\alpha_{22}\triangle_{r,S}\big)\p_t\bigg) \p_t \phi_2 = \frac12\triangle_{r,S}\phi_{2}-g_2+\mu(\Phi_S(t_n)) \phi_2(r, t_n), \quad t_{n}<t<t_{n+1}, \quad n \geq 0, \\
&\Phi_S\left(r, t_{n+1}\right) \triangleq \Phi_S\left(r, t_{n+1}^{+}\right)=\frac{\Phi_S\left(r, t_{n+1}^{-}\right)}{\left\|\Phi_S\left(r, t_{n+1}^{-}\right)\right\|_{L^2(\mathbb{R}^+)}}, \quad n \geq 0,\\
&\Phi_S\left(r, 0\right) = \Phi_S^0(r), \quad \text{with} \quad  \|\Phi_S^0\|^2_{L^2(\mathbb{R}^+)} =1, \quad \dot{\Phi}_S\left(r, 0\right) = \dot{\Phi}_S^0(r),
\end{split}
\end{align}
where
\begin{align}
\begin{split}
g_1\triangleq \big(V_1-\eta-\beta\phi_{2}^2\big)\phi_{1} -\big(H_0+H_1\big)\phi_{2}, \\
g_2\triangleq \big(V_2+\eta-\beta\phi_{1}^2\big)\phi_{2} -\big(H_0+H_1\big)\phi_{1}.
\end{split}
\end{align}

Define $\dot{\phi}_i\triangleq \p_t \phi_{i}, \quad i=1,2$. Then the above system \eqref{binary GFLM with stibilizer} can be rewritten into:
\begin{align} \label{binary GFLM with stibilizer 1}
\begin{split}
&\p_t \phi_{1}=\dot{\phi}_1,\\
&\bigg(\alpha_{01}+\big(\alpha_{11}-\alpha_{21}\triangle_{r,S}\big)\p_t\bigg) \dot{\phi}_1 = \frac12\triangle_{r,S}\phi_{1}-g_1+\mu(\Phi_S(t_n)) \phi_1(r, t_n), \quad t_{n}<t<t_{n+1}, \quad n \geq 0, \\
&\p_t \phi_{2}=\dot{\phi}_2,\\
&\bigg(\alpha_{02}+\big(\alpha_{12}-\alpha_{22}\triangle_{r,S}\big)\p_t\bigg) \dot{\phi}_2 = \frac12\triangle_{r,S}\phi_{2}-g_2+\mu(\Phi_S(t_n))\phi_2(r, t_n), \quad t_{n}<t<t_{n+1}, \quad n \geq 0, \\
&\Phi_S\left(r, t_{n+1}\right) \triangleq \Phi_S\left(r, t_{n+1}^{+}\right)=\frac{\Phi_S\left(r, t_{n+1}^{-}\right)}{\left\|\Phi_S\left(r, t_{n+1}^{-}\right)\right\|_{L^2(\mathbb{R}^+)}}, \quad n \geq 0,\\
&\Phi_S\left(r, 0\right) = \Phi_S^0(r), \quad \text{with} \quad  \|\Phi_S^0\|^2_{L^2(\mathbb{R}^+)} =1,  \quad \dot{\Phi}_S\left(r, 0\right) = \dot{\Phi}_S^0(r).
\end{split}
\end{align}

In  $\mathbb{R}^+$ with confining potential $V_i(r)(i.e., lim_{r\rightarrow +\infty} V_i(r) =+\infty, i=1,2)$, we note that the symmetric and central vortex state solution  decays exponentially fast as $r\rightarrow +\infty$ \cite{Cazenave2003}. Hence, the unbounded domain $\mathbb{R}^+$ can be truncated into a sufficient large bounded interval $U=[0, R]$ when one solves for the steady state solution.
We use the following coordinate transformation $r=\frac{R}2(x+1)$ to transform the equation into one defined in $x\in I\triangleq[-1,1]$ by setting $u_i(x)=\phi_i(\frac{R}2(x+1))$, $v_i(x)=\dot{\phi}_i(\frac{R}2(x+1)), i=1, 2$, and $\tilde{H}_1(x)=H_1(r)$. The transformed equation system is given by
\begin{align}\label{binary scaling BFET}
\begin{split}
&\p_t u_1=v_1, \quad t_{n}<t<t_{n+1}, \quad n \geq 0,\\
&\bigg(\alpha_{01}+\big(\alpha_{11}-\frac{4\alpha_{21}}{R^2} \triangle_{x,S}\big)\p_t\bigg) v_1  = \frac{2}{R^2}\triangle_{x,S} u_1 - \tilde{g}_1 +\mu(\Phi_S(t_n)) u_1(t_n), \quad t_{n}<t<t_{n+1}, \quad n \geq 0,\\
&\p_t u_2=v_2, \quad t_{n}<t<t_{n+1}, \quad n \geq 0,\\
&\bigg(\alpha_{02}+\big(\alpha_{12}-\frac{4\alpha_{22}}{R^2} \triangle_{x,S}\big)\p_t\bigg) v_2  = \frac{2}{R^2}\triangle_{x,S} u_2 - \tilde{g}_2 +\mu(\Phi_S(t_n)) u_2(t_n), \quad t_{n}<t<t_{n+1}, \quad n \geq 0,\\
&\myvec{u}(t_{n+1})\triangleq\frac{\myvec{u}(t^+_{n+1})}{\|\myvec{u}(t^+_{n+1})\|_{L^2(I)}}, \quad n \geq 0,\\
&\myvec{u}(x,0) = \Phi_S^0(r), \quad \text{with} \quad  \|\myvec{u}^0\|^2_{L^2(I)}=1, \quad \myvec{v}(x,0) = \dot{\Phi}_S^0(r),
\end{split}
\end{align}
where
\begin{align}
\begin{split}
\tilde{g}_1\triangleq \big(\tilde{V_1} -\eta-\beta u_2^2\big) u_1 -\big(H_0+\tilde{H}_1\big)u_{2}, \\
\tilde{g}_2\triangleq \big(\tilde{V_2} +\eta-\beta u_1^2\big) u_2 +\big(H_0+\tilde{H}_1\big)u_{1},
\end{split}
\end{align}
and $\|\myvec{u}(t^+_{n+1})\|^2_{L^2(I)}\triangleq \|u_1(t^+_{n+1})\|^2_{L^2(I)} +\|u_2(t^+_{n+1})\|^2_{L^2(I)}$.
The chemical potential $\mu$ at time $t=t_n$ can be rewritten as
\begin{align}
\mu(\Phi_S(t_n))=&\sum_{i=1,2}\pi\int_{-1}^{1}\bigg(\big(\p_x u_i(t_n)\big)^2+\frac{S^2}{(x+1)^2} u_i^{2}(t_n)\bigg)\ (x+1)dx \nonumber\\&+\frac{\pi R^2}{2}\int_{-1}^{1}\bigg(\big(\tilde{V_1} -\eta-\beta u_2^2(t_n)\big) u_1(t_n) +\big(H_0+\tilde{H}_1(t_n)\big)u_{2} (t_n)\bigg)u_{1} (t_n)\ (x+1)dx\nonumber\\&+ \frac{\pi R^2}{2}\int_{-1}^{1}\bigg( \big(\tilde{V_2} +\eta-\beta u_1(t_n)^2\big) u_2(t_n) +\big(H_0+\tilde{H}_1(t_n)\big)u_{1}(t_n)\bigg)u_{2} (t_n)\ (x+1)dx.
\end{align}

Following the  development for ASGF-\uppercase\expandafter{\romannumeral 1\relax}, we obtain the following decoupled discrete schemes.
\begin{algorithm}
Given initial data $u_i^0, v_i^0$, $i=1,2$, compute the spectral coefficients $\vec {\hat{u_i}}^0, \vec{\hat{v_i}}^0$. For $n>0$, compute $\vec{\hat{u_i}}^n, \vec{\hat{v_i}}^n$ via
\begin{align}\label{binary ASGF 1}
\begin{split}
&\frac{\vec{\hat{u}}_1^{*}-\vec{\hat{u}}_1^{n}}{\triangle t}=\vec{\hat{v}}_1^{*} ,\quad t_{n}<t<t_{n+1}, \quad n \geq 0, \\
&\alpha_{01}C\vec{\hat{v}}_1^{*}  + \big(\alpha_{11}C+\frac{4\alpha_{21}}{R^2}(A+S^2 B)\big) \frac{\vec{\hat{v}}_1^{*}-\vec{\hat{v}}_1^{n}}{\triangle t} = \big(-\frac2{R^2}(A+S^2 B)-\alpha_1 C\big)\vec{\hat{u}}_1^{*}+ C\big(\alpha_1\vec{\hat{u}}_1^{n} +\vec{\hat{g}}_1^n+ \mu(u_N^n)\vec{\hat{u}}_1^{n}\big),\\
&\frac{\vec{\hat{u}}_2^{*}-\vec{\hat{u}}_2^{n}}{\triangle t}=\vec{\hat{v}}_2^{*} ,\quad t_{n}<t<t_{n+1}, \quad n \geq 0, \\
&\alpha_{02}C\vec{\hat{v}}_2^{*}  + \big(\alpha_{12}C+\frac{4\alpha_{22}}{R^2}(A+S^2 B)\big) \frac{\vec{\hat{v}}_2^{*}-\vec{\hat{v}}_2^{n}}{\triangle t} = \big(-\frac2{R^2}(A+S^2 B)-\alpha_2 C\big)\vec{\hat{u}}_2^{*}+ C\big(\alpha_2\vec{\hat{u}}_2^{n} +\vec{\hat{g}}_2^n+ \mu(u_N^n)\vec{\hat{u}}_2^{n}\big),\\
&\quad t_{n}<t<t_{n+1}, \quad n \geq 0,\\
&\myvec{u}_N^{n+1}=\frac{\myvec{u}_N^{*}}{\|\myvec{u}_N^{*}\|_{L^2(I)}},\quad n \geq 0, \quad \forall S \geq 0,
\end{split}
\end{align}
where $\alpha_1 = max\{\frac12\big(\tilde{V_1} -\eta-\beta (u^n_2)^2 +|H_0+\tilde{H}_1^n| -\mu(u^n)\big),0\}$  and $\alpha_2 = max\{\frac12\big(\tilde{V_2} +\eta-\beta (u^n_1)^2 +|H_0+\tilde{H}_1^n| -\mu(u^n)\big),0\}$ are chosen stabilization parameters such that the time step can be as large as possible;

until the $L_2$-norm of the residue is less than tolerance $\epsilon$.
\end{algorithm}

\subsection{PPNCG method}

Analogous to the development of numerical schemes for the single-component GP equation, we arrive at the  PPNCG method for the binary GP system.
The update in this method is defined by
\begin{align}\label{alpha_n bcg}
\begin{split}
\tilde{\Phi_{S}}^{n+1} &= \Phi_{S}^{n}  -\alpha_n \ diag\left(P^n, P^n\right) \ \myvec{r}^n,\\
\Phi_S^{n+1}&=\frac{\tilde{\Phi_{S}}^{n+1} }{\|\tilde{\Phi_{S}}^{n+1} \|_{L^2(\mathbb{R}^+)}},
\end{split}
\end{align}
where
\begin{align}
P^n:=\left(\pi\int_0^{\infty}\left(\left((\phi_{1}^n)^{\prime}\right)^{2} +\frac{S^2}{r^2} (\phi_{1}^n)^{2} + \left((\phi_{2}^n)^{\prime}\right)^{2} +\frac{S^2}{r^2} (\phi_{2}^n)^{2}\right)\ rdr - \frac12\triangle_{r,S}\right)^{-1}
\end{align}
is the symmetric positive definite preconditioner, and
\begin{align} \label{residuals}
\myvec{r}^n=\left( \begin{array}{c} r_1^n \\   r_2^n \end{array}\right) \triangleq
\left( \begin{array}{c} -\frac12\triangle_{r,S} \phi_{1}^n+\bigg[V_{1}-\eta-\beta(\phi_{2}^n)^2\bigg]\phi_{1}^n-\big(H_0+H_1^n\big)\phi_{2}^n
- \mu^n \phi_{1}^n\\  -\frac12\triangle_{r,S}\phi_{2}^n+\bigg[V_{2}+\eta-\beta(\phi_{1}^n)^2\bigg]\phi_{2}^n-\big(H_0+H_1^n\big)\phi_{1} ^n- \mu^n \phi_{2}^n \end{array}\right)
\end{align}
is the residue.

We reformulate this formula as follows
\begin{align}\label{theta_n bcg}
\Phi_S^{n+1} = cos(\theta_n) \Phi_S^{n} + sin(\theta_n) \frac{\myvec{p}^n}{\|\myvec{p}^n\|_{L^2(\mathbb{R}^+)}},
\end{align}
with
\begin{align}
\myvec{p}^n = \myvec{d}^n -\ \langle \myvec{d}^n, \Phi_{S}^{n} \rangle \Phi_S^{n},
\end{align}
 where $\myvec{d}^n= - diag\left(P^n, P^n\right) \myvec{r}^n+\beta_n \myvec{d}^{n-1}$ is the search direction.
Expanding $\Phi_S^{n+1}$ up to second-order in $\theta_n$, we obtain
\begin{align}
\Phi_S^{n+1}=\left(1-\frac{\theta_{n}^{2}}{2}\right) \Phi_S^{n}+\theta_{n} \frac{\myvec{p}^n}{\|\myvec{p}^n\|_{L^2(\mathbb{R}^+)}} +\mathcal{O}\left(\theta_{n}^{3}\right)
\end{align}
and
\begin{align}\label{E(Phi_n1)}
E\left(\Phi_S^{n+1}\right)=E\left(\Phi_S^{n}\right)+\frac{2\theta_{n}}{\|\myvec{p}^n\|} \left\langle \myvec{r}^n, \myvec{p}^n\right\rangle +\frac{\theta_{n}^{2}}{\|\myvec{p}^n\|^2}\left[\frac{1}{2}  \nabla^2 E\left(\Phi_S^{n} \right)[\myvec{p}^n, \myvec{p}^n]-\mu^{n}\|\myvec{p}^n\|^2\right]+\mathcal{O}\left(\theta_{n}^{3}\right).
\end{align}

Minimizing \eqref{E(Phi_n1)} with respect to $\theta_{n}$ yields
\begin{align}
\theta_{n}= - \frac{\|\myvec{p}^n\| \cdot\left\langle \myvec{r}^n, \myvec{p}^n\right\rangle}{\frac{1}{2} \nabla^2 E\left(\Phi_S^{n} \right)[\myvec{p}^n, \myvec{p}^n] -\mu^{n}\|\myvec{p}^n\|^2}.
\end{align}

 The algorithm is a straight-forward extension of the PPNCG method for the single-component GP equation model, which we will not repeat here.

\subsection{Numerical results of the binary GP model}

We first benchmark the ASGF method against the GFLM method for the binary GP model and then compare the ASGF method with the PPNCG method. Finally, we present some steady state solutions computed using the PPNCG for some selected model parameters.

\noindent{\bf Example 3.1.}
 We compare the $ASGF-\uppercase\expandafter{\romannumeral 1\relax}$ method with the GFLM method (i.e., using $\alpha_0=1, \alpha_1=\alpha_2=0$ in the ASGF-\uppercase\expandafter{\romannumeral 1\relax} scheme) when computing the central vortex state solution.
The initial datum is chosen as  $\Phi_S^0(r) = (\sqrt{\alpha}\Phi^0(r), \sqrt{1-\alpha}\Phi^0(r))^T$, where $\Phi^0(r) = \frac 1{\sqrt{\pi S!}} r^S e^{-r^2/2}$ and $\dot{\Phi_S}^0 = 0$. The stopping criterion in time marching is that the $L^2$-norm of residue of the Euler–Lagrange equation \eqref{binary lagrange equation} is less than tolerance $\varepsilon = 10^{-10}$.
We take  $\eta=10$, $\gamma=\pi$, $\alpha=0.5$, $H_0=5$, $U=[0,16]$ and $V_1(r)=V_2(r)=\frac {r^2}2+25sin^2(\frac{\pi r}{4})$ with the dimension of the discrete Legendre space $N_l = 200$. The performance of the schemes with $S=3$, $\beta=60$ and with $S=7$, $\beta=100$, using different time steps $\tau$, is shown in Table \ref{Tab31} and Table \ref{Tab32}, respectively, where $E_c$ and $\mu_c$ are the energy and chemical potential of the central vortex state solution,  and \#$iter$ is the number of iterations.

From the numerical results, we observe that the computational time of the GFLM method at  $\alpha_0=1, \alpha_1=\alpha_2=0$ is much longer than that of the ASGF method at the same time step. Hence, the stabilizer $\alpha_0+\big(\alpha_1-\alpha_2 \triangle_{r,S}\big)\p_t$ indeed speeds up the computation to a quite large extent.  The results show that the ASGF method is more efficient than the GFLM method.

\begin{table}
	\centering
	\begin{tabular}{|c|ccccccc|}
		\hline
		$\tau$ &$\alpha_0$& $\alpha_1$&$\alpha_2$&$\mathrm{CPU}(\mathrm{s})$ & $ E_{c}$ & $\mu_{c} $& \#$iter$ \\
		\hline
			   &1&0&0&15.61&-0.5052747150&-0.5983534336&804\\
				0.01&1E-6&1E-4&1E-4&5.89&-0.5052747150&-0.5983534336 &311\\
				&0&0.01&0&5.37&-0.5052747150&-0.5983534336&276 \\
				\hline
				&1&0&0&7.05&-0.5052747150&-0.5983534336&364\\
				0.1&1E-5&1&0&5.27&-0.5052747150&-0.5983534336& 276\\
				&1E-5&1.25&0.035&4.86&-0.5052747150&-0.5983534336 &247\\
				\hline
				&1&0&0&6.22&-0.5052747150&-0.5983534336&320\\
				1&1E-3&200&3&4.45&-0.5052747150&-0.5983534336&234 \\
				&1E-3&150&3&4.18&-0.5052747150&-0.5983534336 &223\\
				\hline
				\end{tabular}
		\caption{Performance of ASGF-\uppercase\expandafter{\romannumeral 1\relax} scheme when computing the  central vortex state solution with respect to different time step $\tau$ at $S=3, \beta=60$.}
		\label{Tab31}
		\begin{tabular}{|c|ccccccc|}
			\hline
			$\tau$ &$\alpha_0$& $\alpha_1$&$\alpha_2$&$\mathrm{CPU}(\mathrm{s})$ & $ E_{c}$ & $\mu_{c} $& \#$iter$ \\
			\hline
		&1&0&0&14.10&0.7572177467&0.5477025939&798\\
		0.01&1E-7&1.5E-5&1E-4&5.65&0.7572177467&0.5477025939&312\\
		&1E-6&1E-4&5E-4&5.24&0.7572177467&0.5477025939&296 \\
		\hline
		&1&0&0&6.80&0.7572177467&0.5477025939&364\\
		0.1&5E-5&1.5E-4&0.05&5.29&0.7572177467&0.5477025939& 302\\
		&1E-4&1.5E-3&0.05&5.22&0.7572177467&0.5477025939&297\\
		\hline
		&1&0&0&6.21&0.7572177467&0.5477025939&320\\
		1&1E-3&1&5&5.18&0.7572177467&0.5477025939&298 \\
		&8E-3&1.25&5&5.13&0.7572177467&0.5477025939 &296\\
		\hline
	\end{tabular}
	\caption{Performance of ASGF-\uppercase\expandafter{\romannumeral 1\relax} scheme when computing the  central vortex state solution with respect to different time step $\tau$ at  $S=7, \beta=100$. }
	\label{Tab32}
\end{table}

\begin{figure}
		\centering
		\includegraphics[width=10cm]{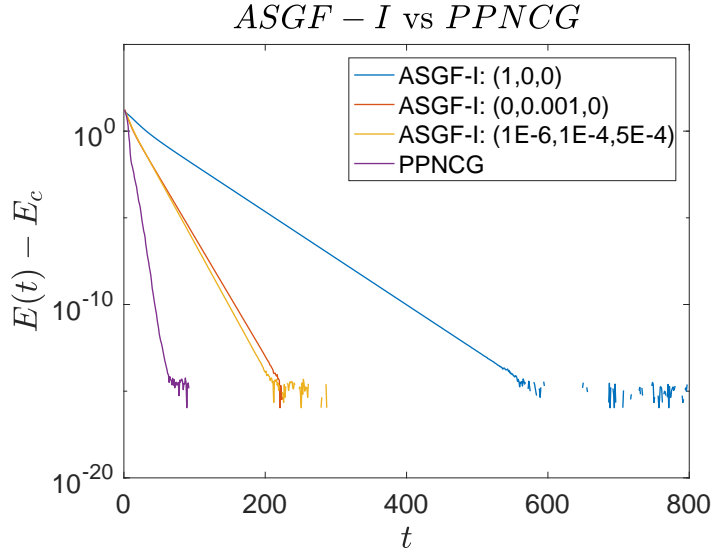}
	\caption{Time evolution of relative energy $E(\Phi_S(·, t))-E_c$ (in logarithmic scale) by different numerical schemes when computing the vortex steady state solution. The triplet, for example $(1,0,0)$ in the inset, represents $\alpha_{0i}=1, \alpha_{1i}=0, \alpha_{2i}=0 \ (i=1, 2)$. The PPNCG method outperforms all others.}
	\label{binaryGFvsPPNCG}
\end{figure}

\noindent{\bf Example 3.2.}
In this example, we compare the total computational time and the number of iterations with respect to different values of $S$ and $\beta$ between ASGF-\uppercase\expandafter{\romannumeral 1\relax} and the PPNCG method when computing the symmetric state and central vortex state solution.
The initial datum is chosen as  $\Phi_S^0(r) = (\sqrt{\alpha}\Phi^0(r), \sqrt{1-\alpha}\Phi^0(r))^T$, where $\Phi^0(r) = \frac 1{\sqrt{\pi S!}} r^S e^{-r^2/2}$ and the initial velocity $\dot{\phi_S}^0 = 100 \times \phi_S^0(r)$, $\eta=50$, $\gamma=\pi$, $\alpha=0.5$, $H_0=50$, $U=[0,16]$  with the dimension of the discrete Legendre space $N_l = 160$ and $V_1(r)=V_2(r)=\frac {r^2}2+25sin^2(\frac{\pi r}{4})$, time step $\tau=1$, $\alpha_{0i}=0.001, \alpha_{1i}=1, \alpha_{2i}=5 \ (i=1, 2)$  for ASGF-\uppercase\expandafter{\romannumeral 1\relax}. The stopping criterion is that the $L^2$-norm of residue of the Euler–Lagrange equation \eqref{binary lagrange equation} is less than  tolerance $\varepsilon = 5\times 10^{-10}$.
We take $S=0, 3, 12, 15$, respectively. The performance of the two  methods is shown in Tab. \ref{Tab33}. The time evolution of the relative energy $E(\Phi_S(·, t))-E_c$ (in logarithmic scale) by different numerical schemes when computing the vortex steady state solution is depicted in Fig. \ref{binaryGFvsPPNCG}. From Tab. \ref{Tab33} and  Fig. \ref{binaryGFvsPPNCG}, we see that the ASGF method outperforms the GFLM method while the PPNCG method is much better than the ASGF-\uppercase\expandafter{\romannumeral 1\relax} method.

\begin{table}
	\centering
	\begin{tabular}{|cc|ccc|ccc|ccc|ccc|}
		\hline
		\multicolumn{2}{|c|}{S}&&0&& &5&&& 10&&&15& \\
		\hline
		\multicolumn{2}{|c|}{$\beta$} &0&5&12&0&100&220&0&200&450&0&300&650\\
		\hline
		\multirow{2}{*}{ASGF-\uppercase\expandafter{\romannumeral 1\relax}} &CPU(s)&3.62&2.99&3.41&4.70&4.94&5.80&4.62&5.77&9.03&4.22&4.10&15.43\\
		&\#$iter$&357&295&342&443&490&570&455&567&904&418&609&1534\\
		\hline
		\multirow{2}{*}{PPNCG} &CPU(s)&1.98&2.03&2.32&2.17&2.36&2.58&1.77&2.05&2.36&1.43&1.60&2.13\\
		&\#$iter$&133&138&157&144&155&163&114&134&160&94&104&143\\
		\hline
	\end{tabular}
	\caption{ Performance comparison between the ASGF-\uppercase\expandafter{\romannumeral 1\relax} method and the PPNCG method in the total computational time and the number of iterations with respect to selected values of $S$ and $\beta$.}
	\label{Tab33}
\end{table}

\begin{figure}[htb!]
	\subfigure[]{
		\begin{minipage}[t]{0.45\linewidth}
			\centering     
			\includegraphics[width=\textwidth]{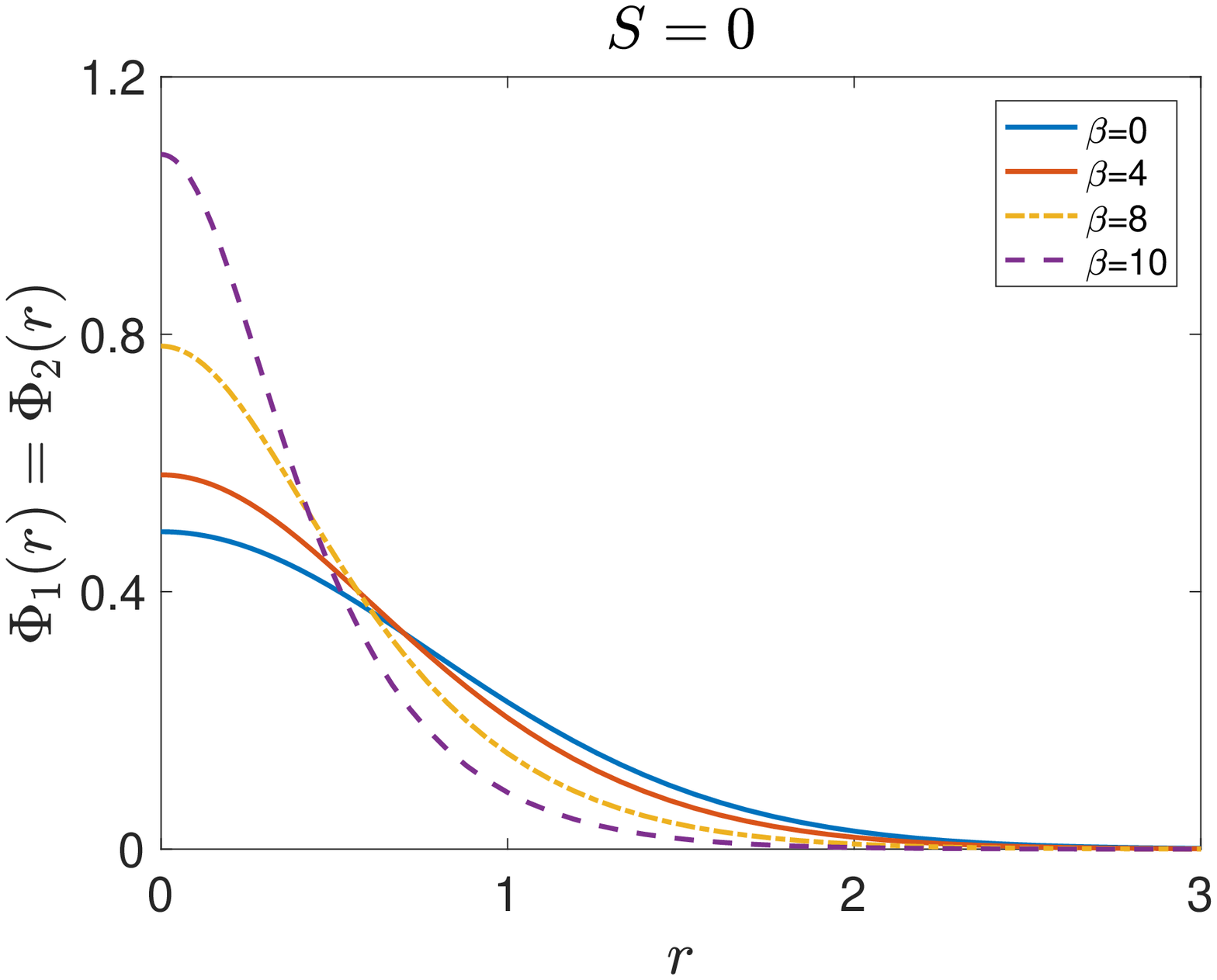}
			\label{fig2a}
	\end{minipage} }
	\hfill
	\subfigure[]{
		\begin{minipage}[t]{0.45\linewidth}
			\centering
			\includegraphics[width=\textwidth]{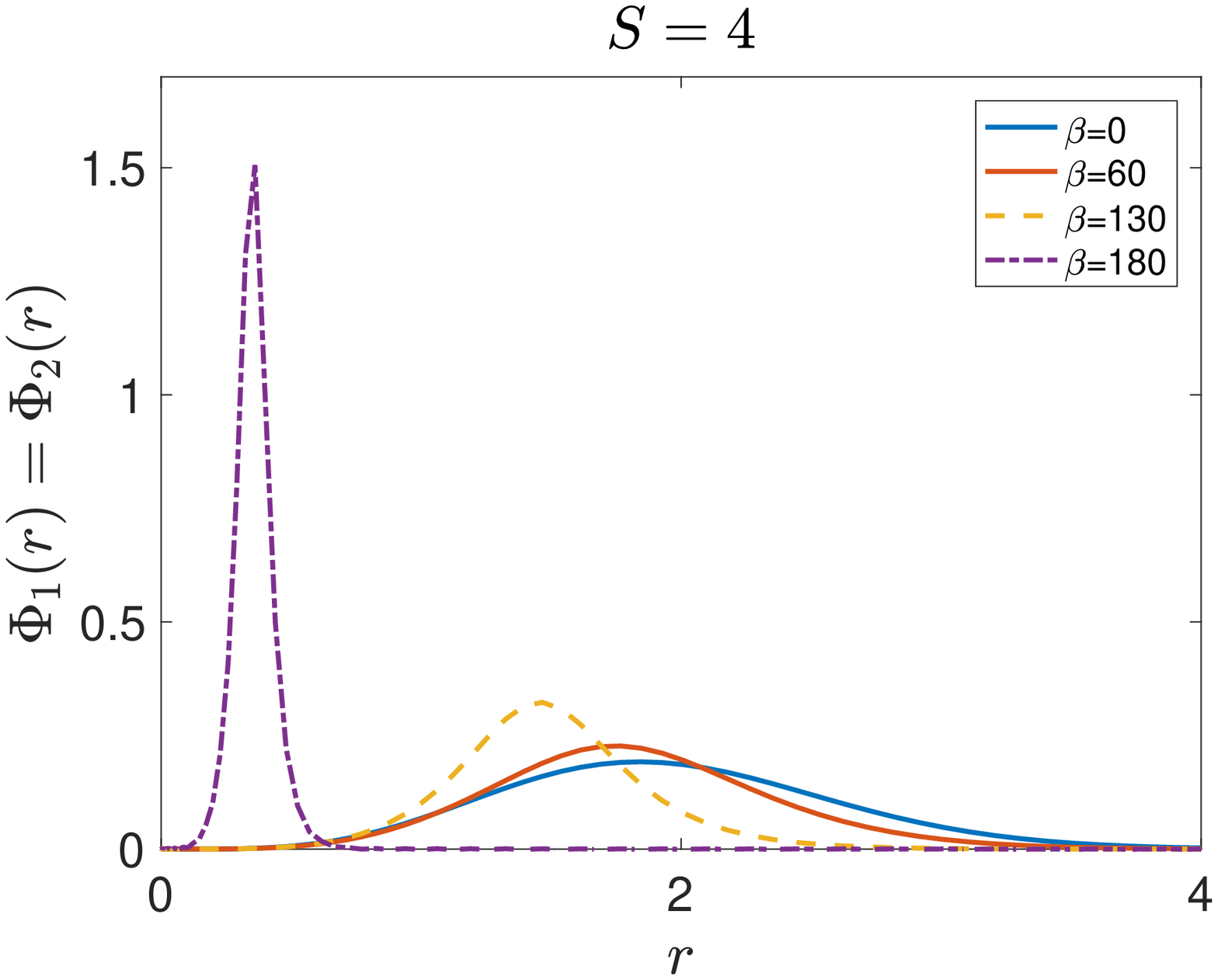}
			\label{fig2b}
	\end{minipage}}
	\subfigure[]{
		\begin{minipage}[t]{0.45\linewidth}
			\centering     
			\includegraphics[width=\textwidth]{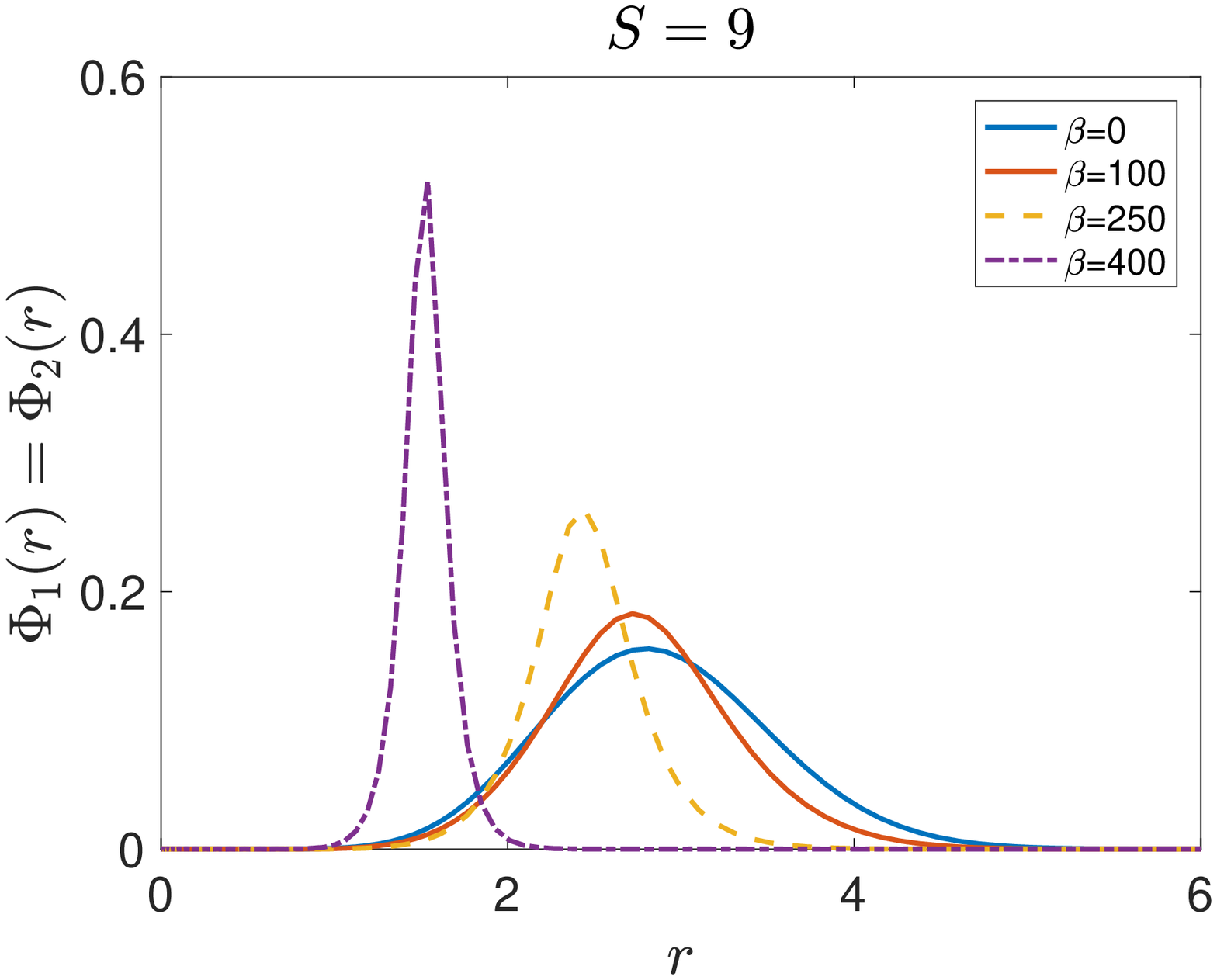}
			\label{fig2c}
	\end{minipage} }
	\hfill
	\subfigure[]{
		\begin{minipage}[t]{0.45\linewidth}
			\centering
			\includegraphics[width=\textwidth]{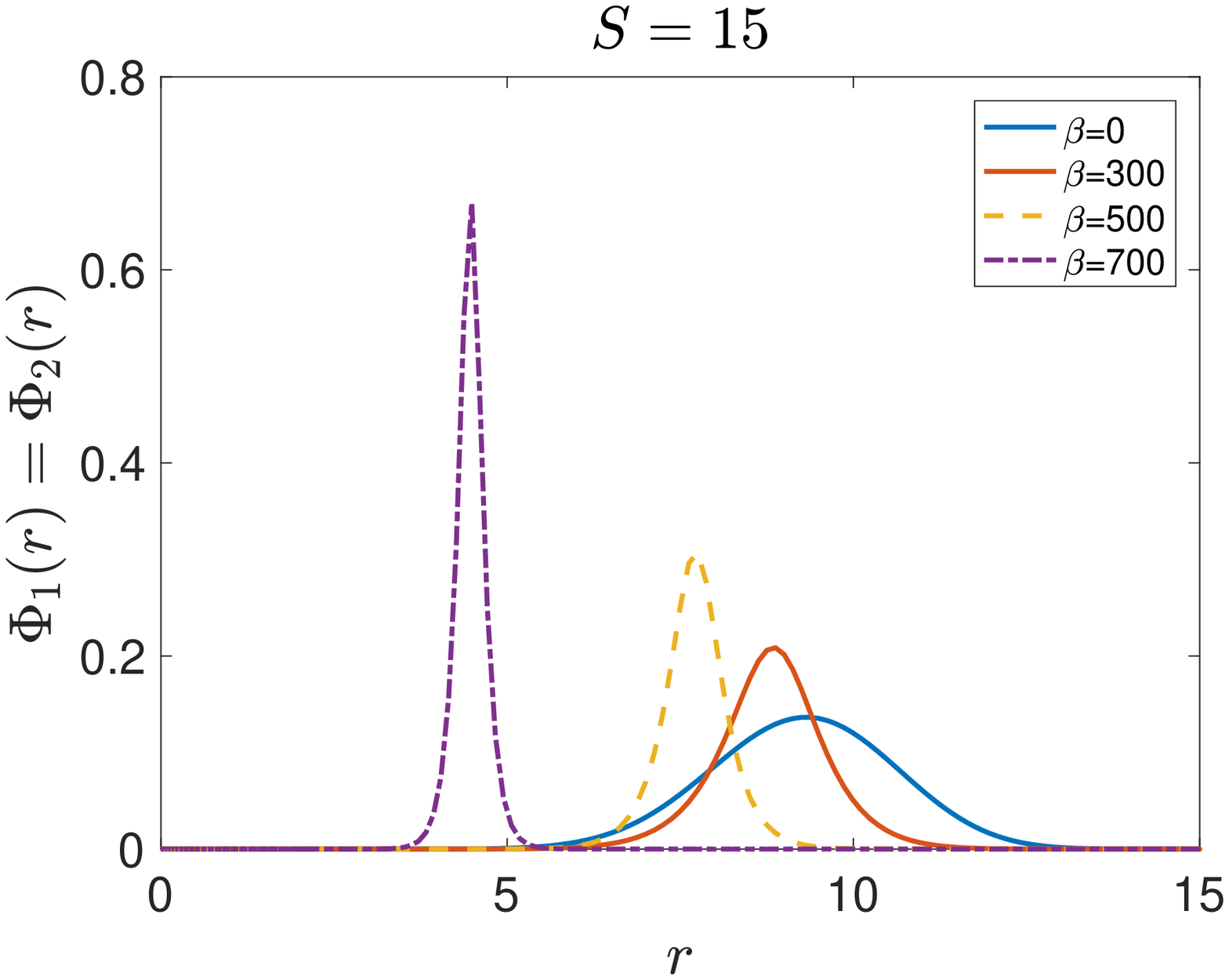}
			\label{fig2d}
	\end{minipage}}
	
	\caption{The annular wave function. The radially dependent wave function is depicted with respect to winding number S=0,4,9,15, respectively. As $\beta$ increases, the radius of the ring expands. }
	\label{ground radial profiles vs S}
\end{figure}

\begin{figure}
	\begin{minipage}[t]{0.45\linewidth}
		\centering     
		\includegraphics[width=\textwidth]{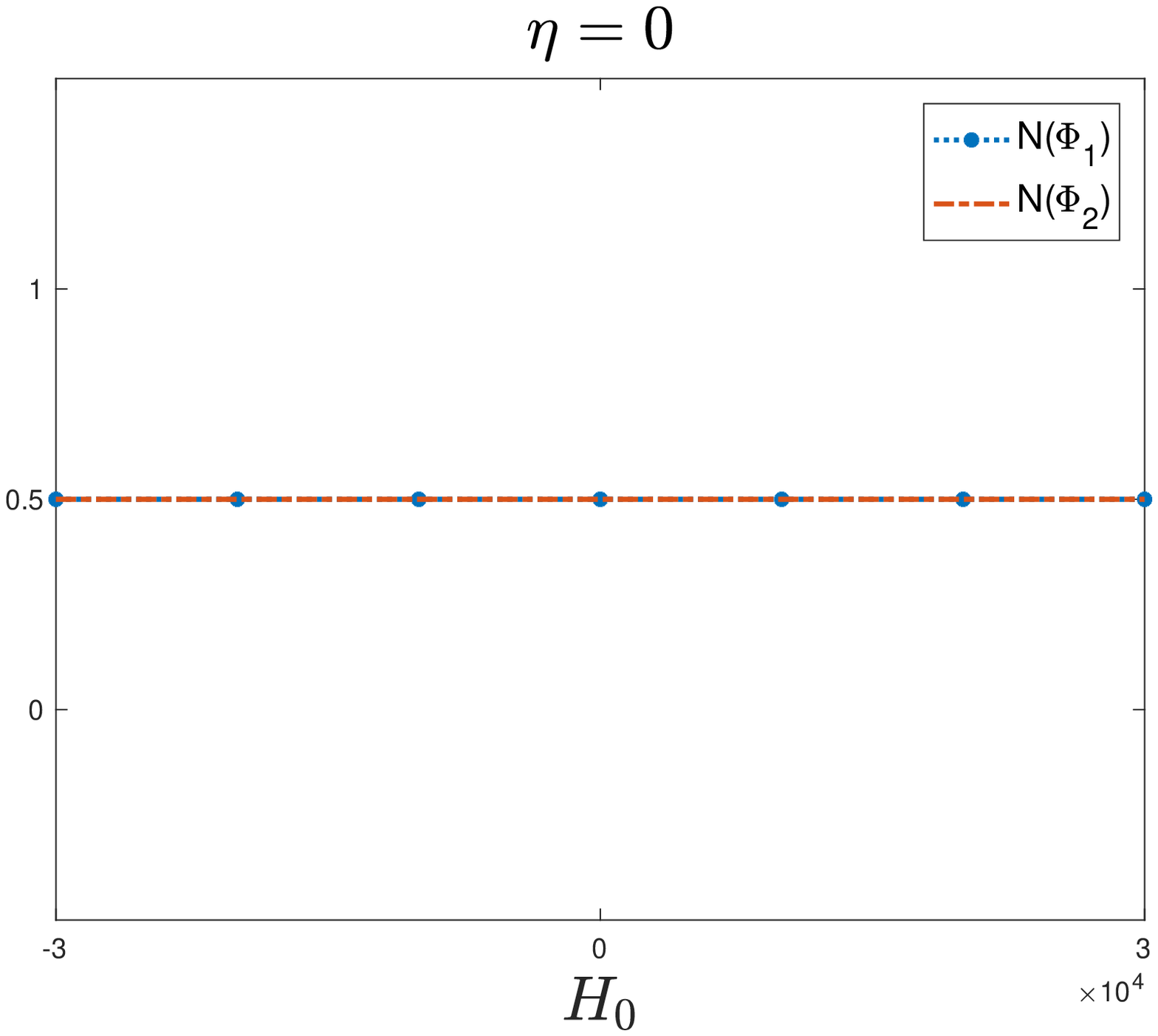}
	\end{minipage} 
	\hfill
	\begin{minipage}[t]{0.45\linewidth}
		\centering
		\includegraphics[width=\textwidth]{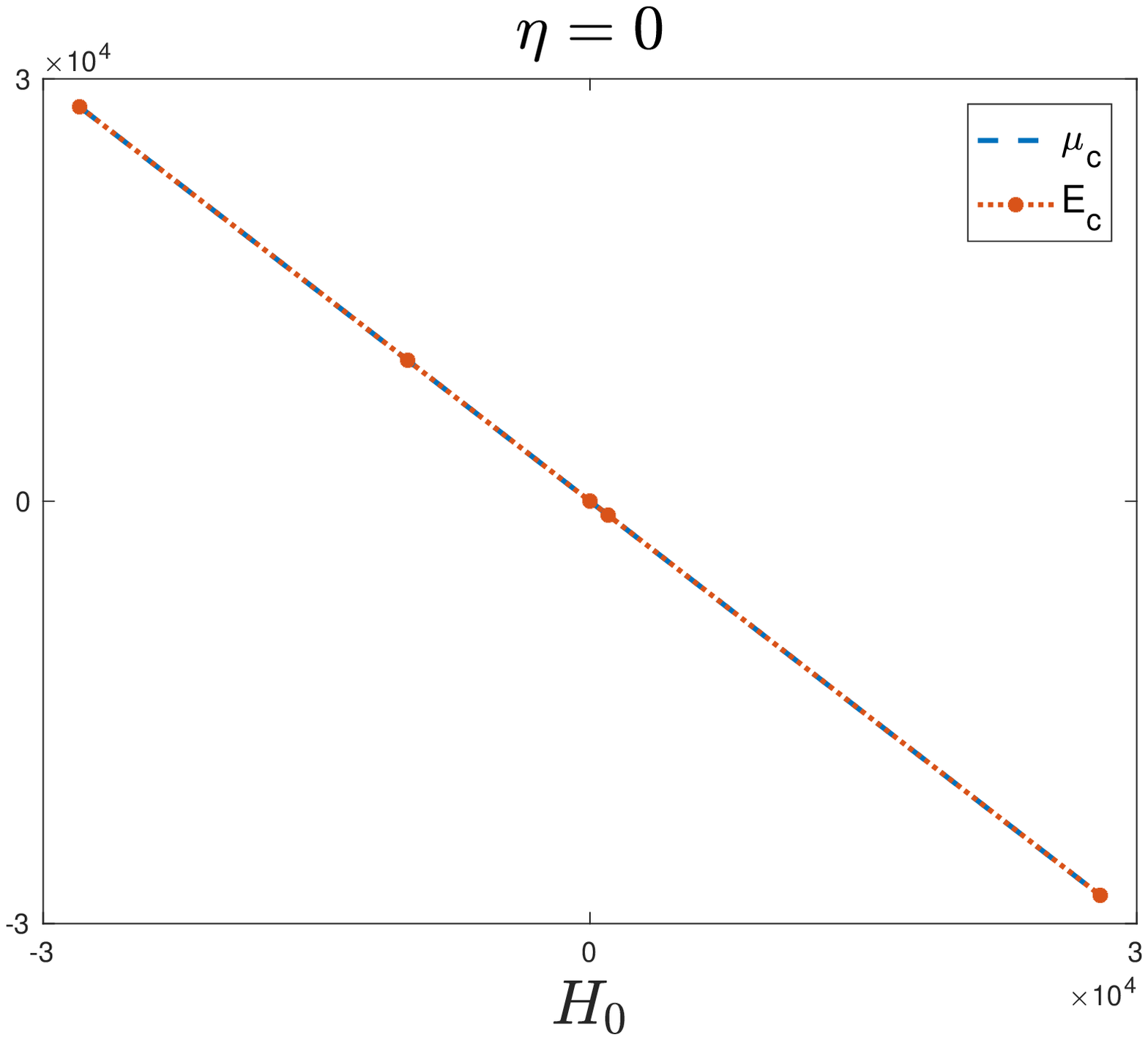}
	\end{minipage}
	\begin{minipage}[t]{0.45\linewidth}
		\centering     
		\includegraphics[width=\textwidth]{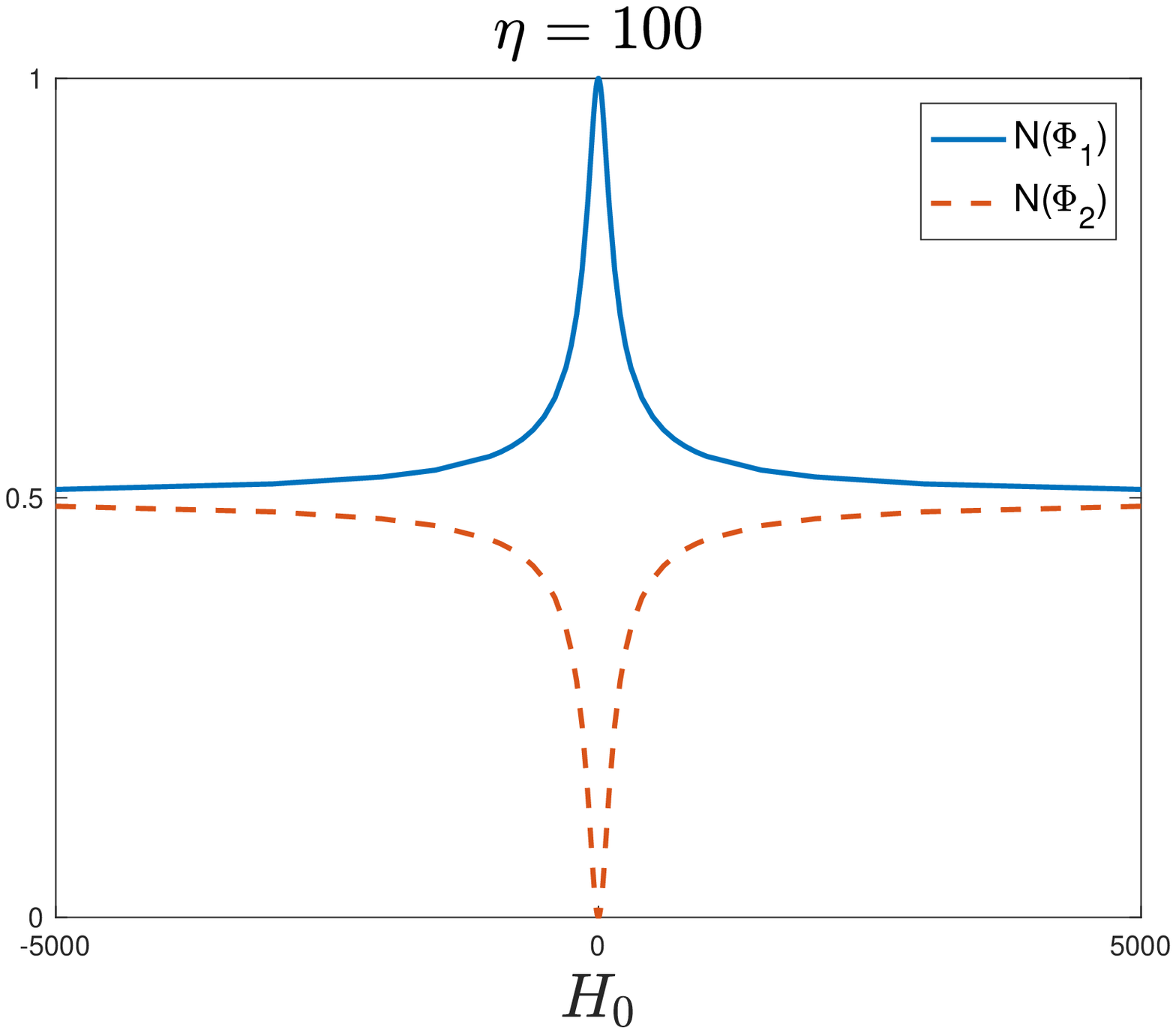}
	\end{minipage} 
	\hfill
	\begin{minipage}[t]{0.45\linewidth}
		\centering
		\includegraphics[width=\textwidth]{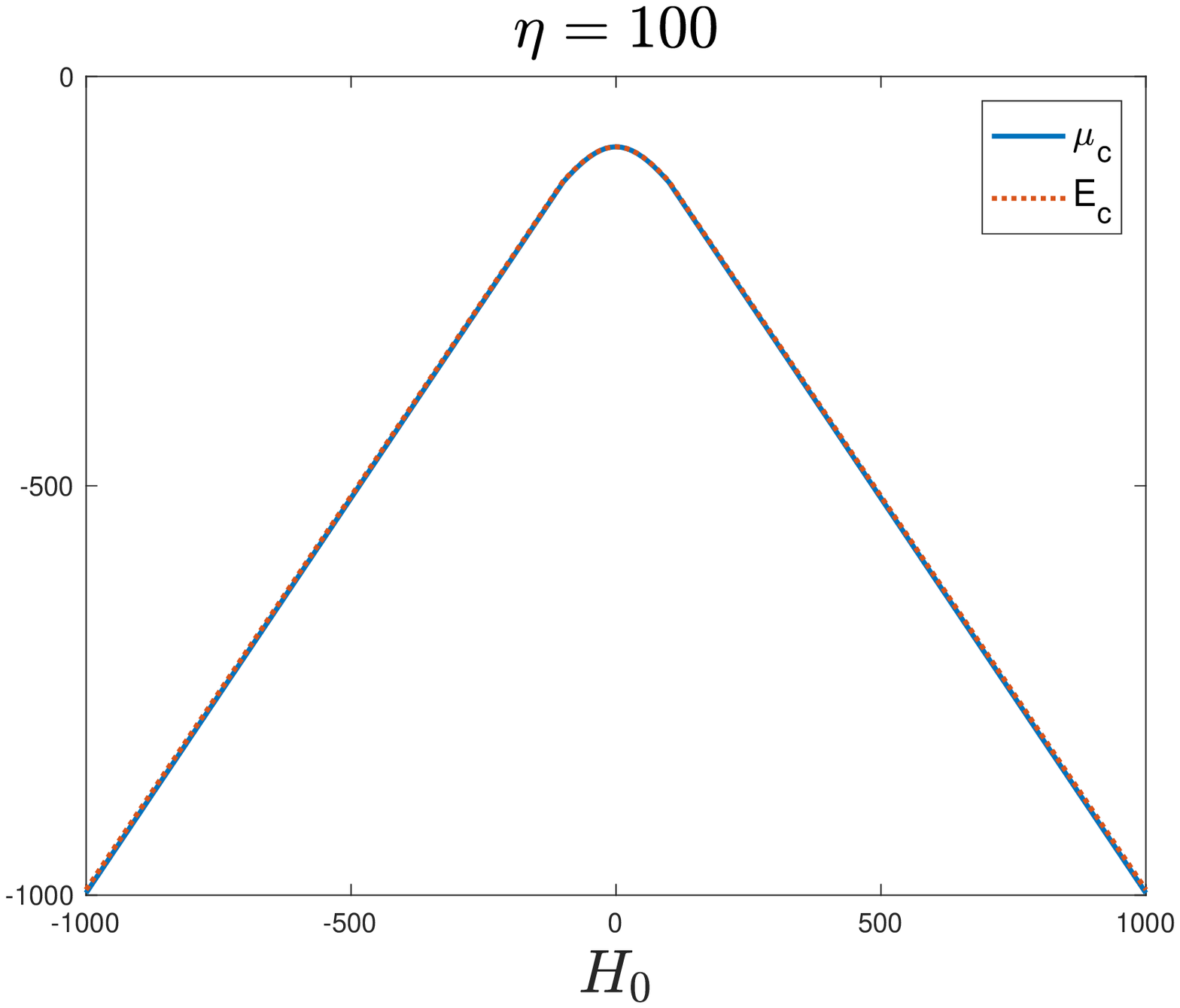}
	\end{minipage}
	\caption{Mass of each component $N(\phi_j) = \|\phi_j\|^2, j =1, 2$, energy $E_c$ and chemical potential $\mu_c$ in vortex steady states.  Variation of mass, energy and chemical potential with respect to $H_0$ at $S=13$ and $\beta=300$ when $\eta= 0$ (Upper) and  when $\eta= 100$ (Lower). The role of the microwave detuning number $\eta$ dominate the mass difference between two states.}
	\label{Fig31}
\end{figure}

\begin{figure}
	\begin{minipage}[t]{0.45\linewidth}
		\centering     
		\includegraphics[width=\textwidth]{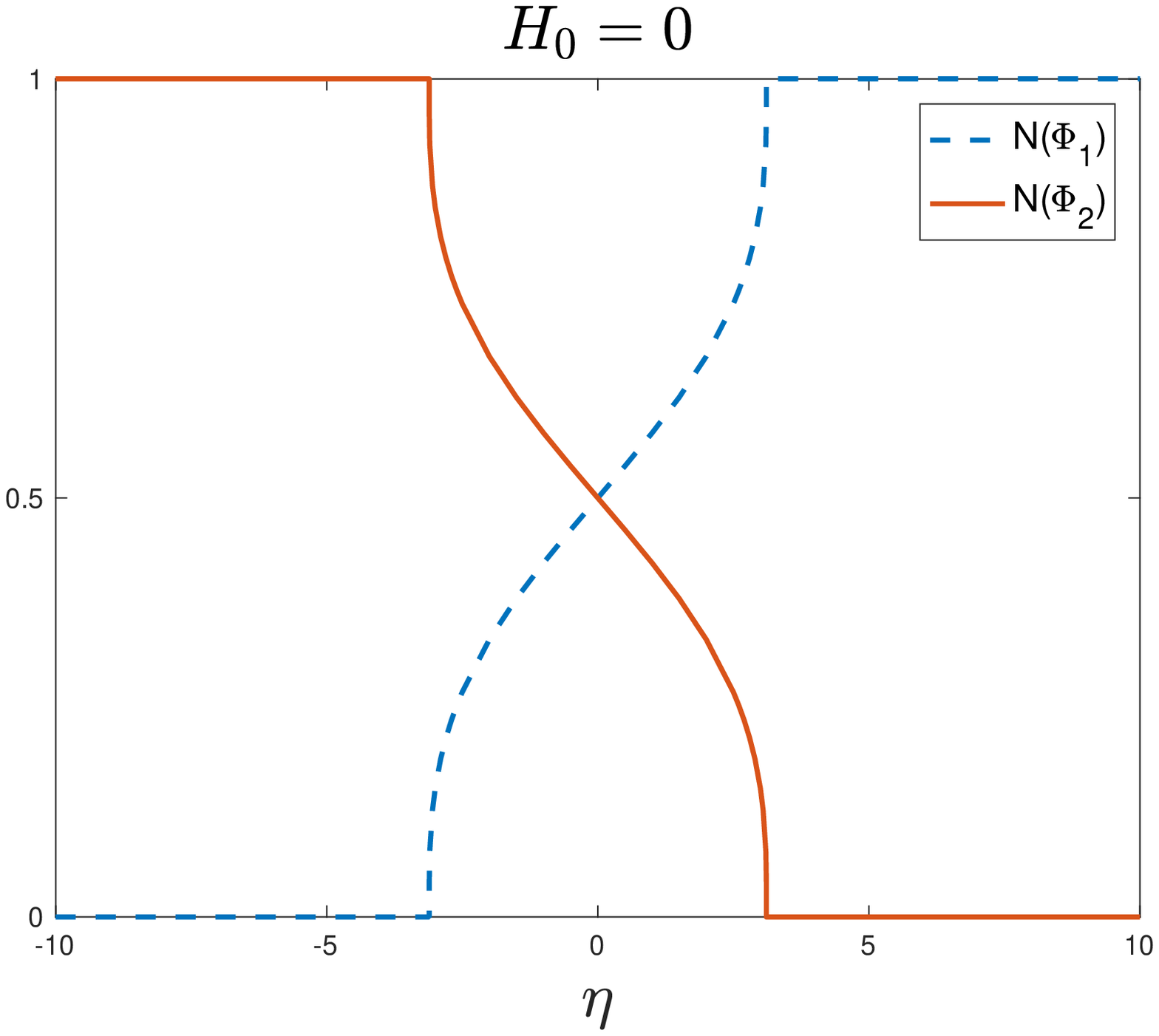}
	\end{minipage} 
	\hfill
	\begin{minipage}[t]{0.45\linewidth}
		\centering
		\includegraphics[width=\textwidth]{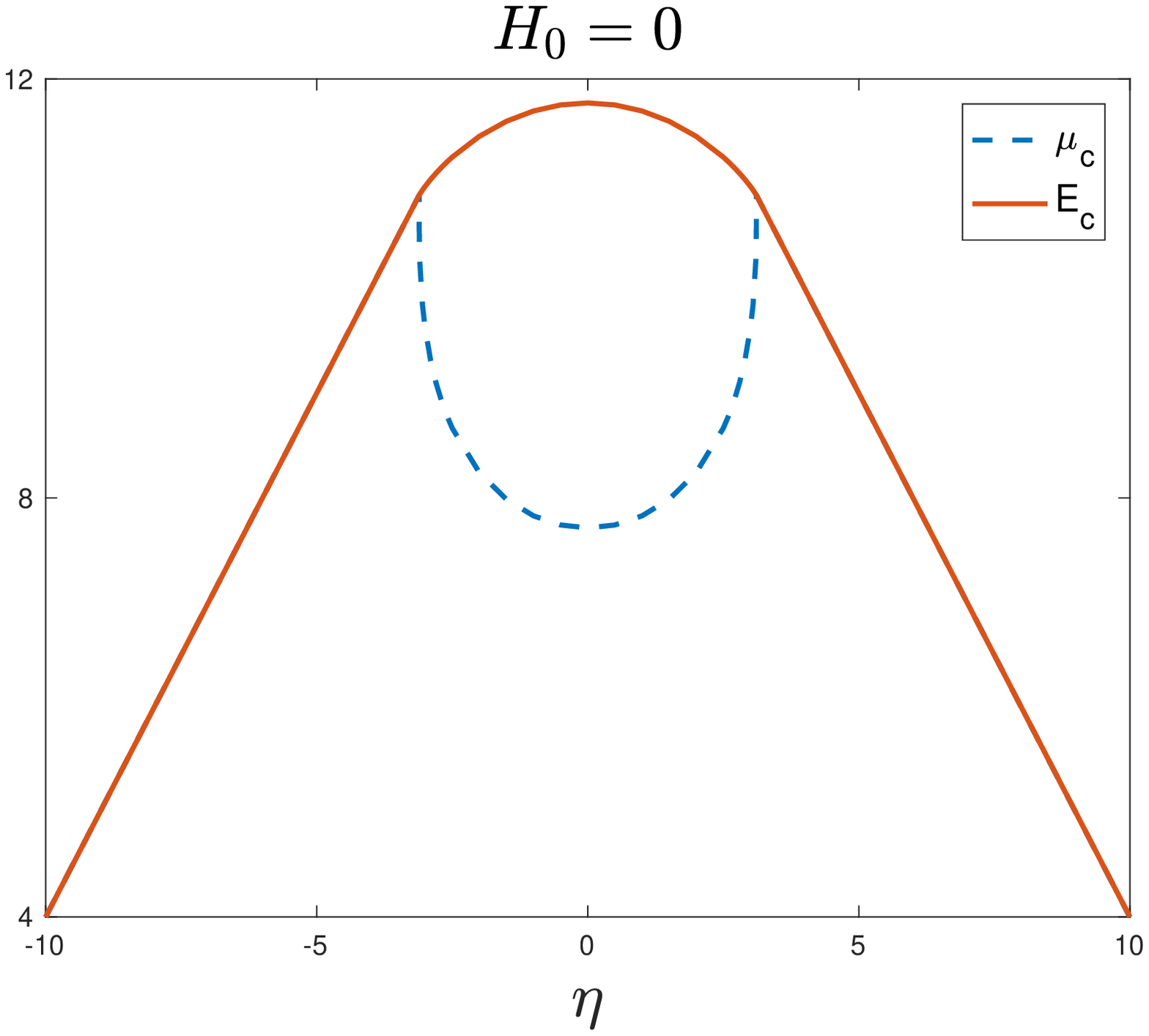}
	\end{minipage}
	\begin{minipage}[t]{0.45\linewidth}
		\centering     
		\includegraphics[width=\textwidth]{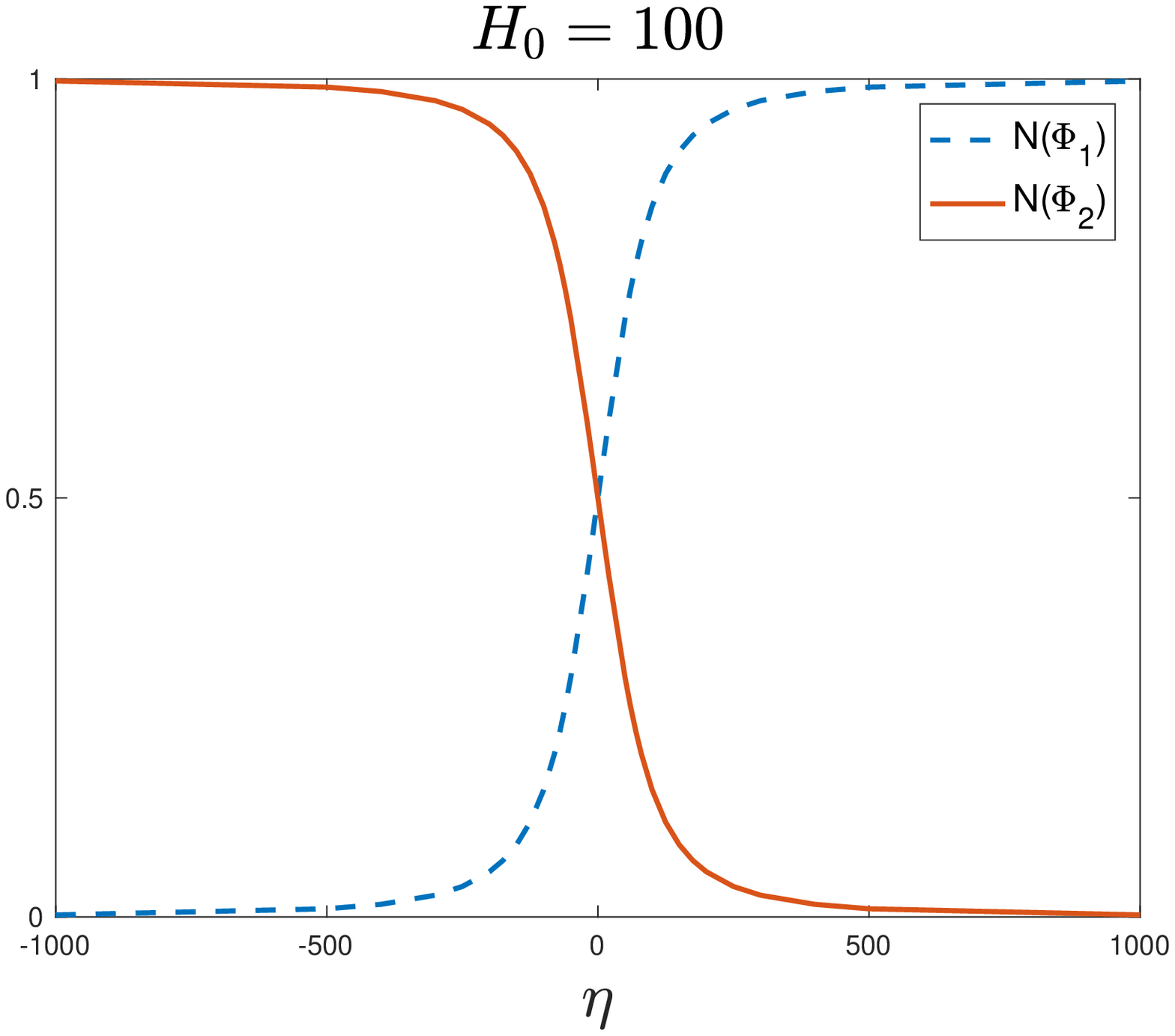}
	\end{minipage} 
	\hfill
	\begin{minipage}[t]{0.45\linewidth}
		\centering
		\includegraphics[width=\textwidth]{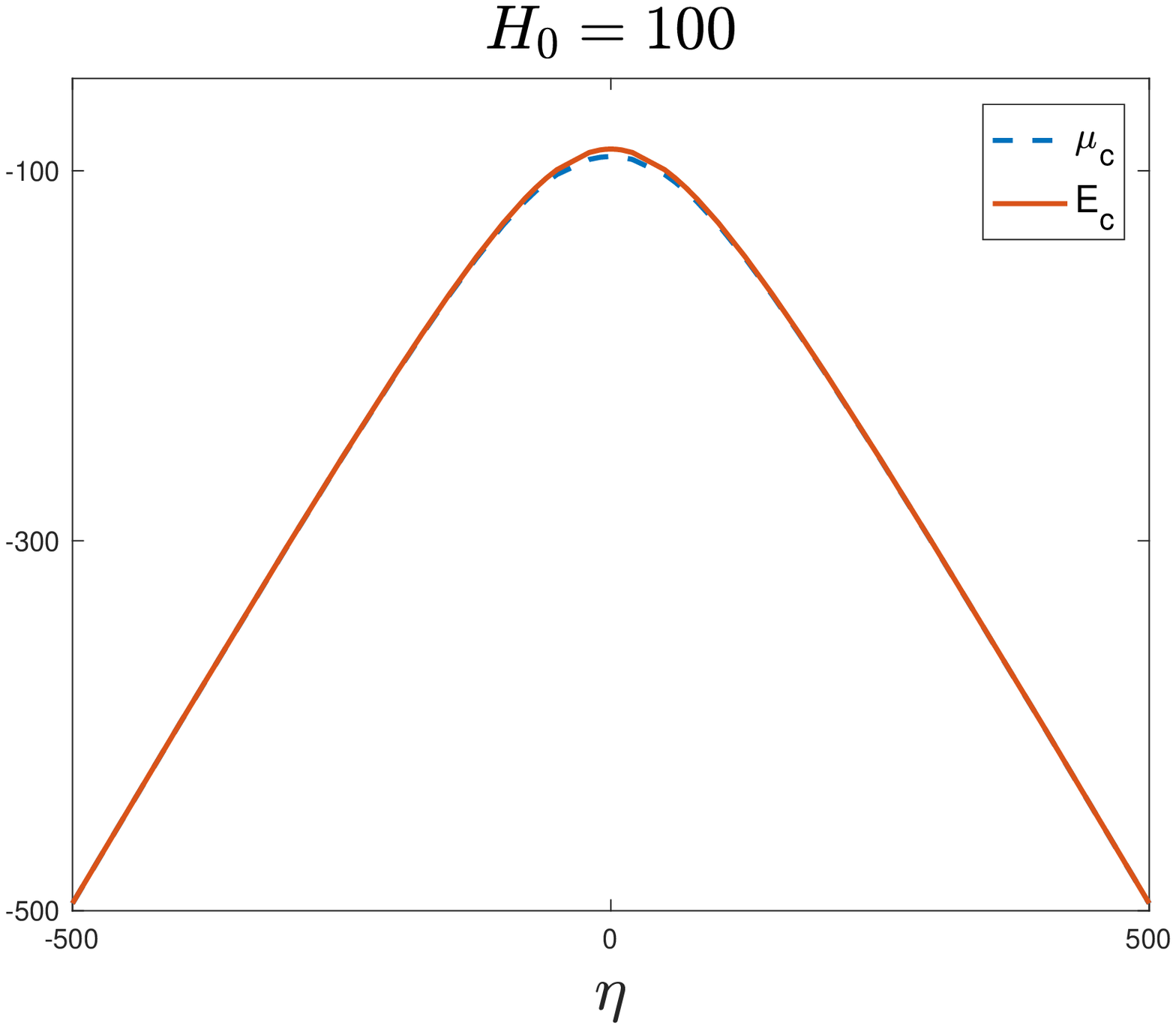}
	\end{minipage}
	\caption{Mass of each component $N(\phi_j) = \|\phi_j\|^2, j =1, 2$, energy $E_c$ and chemical potential $\mu_c$ in vortex steady states.  Variation of mass, energy and chemical potential with respect to $\eta$ at $S=13$ and $\beta=300$ when $H_0 = 0$ (Upper) and  when $H_0 = 100$ (Lower). The role of the background magnetic field $H_0$ is to smooth the mass and chemical potential curve at the transition values so that the other component is not absolutely absent at large $\eta$ and small $\beta$.  }
	\label{Fig32}
\end{figure}

\noindent{\bf Example 3.3.}
Finally, We apply the PPNCG method to obtain some numerical solutions of the binary GP model.
We choose initial datum $\Phi_S^0(r) = (\sqrt{\alpha}\Phi^0(r), \sqrt{1-\alpha}\Phi^0(r))^T$, where $\Phi^0(r) = \frac 1{\sqrt{\pi S!}} r^S e^{-r^2/2}$,  $\eta= 0$, $H_0= 10$, $\alpha = 0.8$, $\gamma = 5\pi$,  and $V_1(r)=V_2(r) = \frac {\gamma}{8\pi}r^2$. We use the same stopping criterion as alluded to early. In Fig. \ref{ground radial profiles vs S}, we  depict  numerical results of  symmetric states $\Phi_{S}^s$  and  vortex steady states $\Phi_{S}^c$ with some selected  values of $S$ and $\beta$.  In Fig. \ref{Fig31}  and \ref{Fig32}, we show changes of mass in each component ($N(\phi_j) = \|\phi_j\|^2, j =1, 2$), energy $E_c:=E(\Phi_S^c)$, and chemical potential $\mu_c :=\mu(\Phi_S^c)$of the   vortex steady states with respect to different  microwave detuning parameter  $\eta$ and background magnetic field $H_0$.

From Fig. \ref{ground radial profiles vs S}, we observe that (i). as the strength of  interaction $\beta$ increases, the radius of vortex annuli in the ground state decreases;  (ii). as winding number $S$ increases, the concentrated (peak) density increases as well.  The results show that $\eta$ is nearly proportional to the mass difference between two states. Whenever $\eta = 0$,  $N(\phi_1) \equiv N(\phi_2)$. When $|\eta|>>1$, one of the two components of vortex steady states dominates and becomes the only possible state in the limit. From Fig. \ref{Fig32}, we observe that there is only one state when $H_0 = 0$ for some small $\beta$ or large $|\eta|$. However, there are two states when $H_0 \neq 0$. $S=0$ is exactly the one that makes the energy the smallest in comparison with $S\neq 0$, which implies the ground  state is spherically symmetric in the general case. 

\section{Conclusion}

We have developed two constrained minimization algorithms for  two steady GP equations coupled with magnetic field, equivalent to two nonlocal GP systems, based on the normalized gradient flow model and the perturbed, projected conjugate gradient approach, respectively. These methods are firstly presented using the single component GP model, and later extended to the binary GP model. Detailed comparisons among the new algorithms and the existing GFLM method when computing  symmetric and central vortex state solutions of the GP models are conducted
in 2D space. The comparative study shows that the ASGF method is significantly better than the GFLM method, while the PPNCG scheme outperforms the ASGF scheme in all the cases investigated. These new methods can be readily extended to other GP equation systems with different external potentials or without the magnetic field coupling, adding additional, efficient  computational tools for GP systems.

\section*{Acknowledgements}
 Di Wang's research is partially supported by 
 NSAF-U1930402 awarded to CSRC.
Qi Wang's research is partially supported by NSF awards  OIA-1655740 and a GEAR award from SC EPSCoR/IDeA Program.

\bibliographystyle{plain}
\bibliography{Ref}

\end{document}